%% file: main.tex
\documentclass[11pt]{article}

\usepackage[T1]{fontenc}
\usepackage{textcomp}
\usepackage[main=english, italian]{babel}
\usepackage[a4paper,margin=3cm]{geometry}
\usepackage{graphicx}
\usepackage{indentfirst}
\usepackage{microtype}
\usepackage{quoting}
\usepackage{csquotes}
\usepackage{amsfonts}
\usepackage{amsmath}
\usepackage{caption}
\usepackage{pdfpages}
\usepackage{siunitx} 
\usepackage{graphicx}
\usepackage{mathrsfs}
\usepackage{comment}
\usepackage{mathtools}
\usepackage{float}

\usepackage[backend=biber]{biblatex}
\usepackage{hyperref}
\usepackage{hyperref}
\usepackage{authblk}
\usepackage{physics}

\graphicspath{{images/}}

\makeindex

\usepackage[absolute]{textpos}
\usepackage{xcolor}

\newcommand{\keywords}[1]{%
    \vspace{1em}
    \noindent\textbf{Keywords:} #1
}

\author[1]{Elena Scibona}
\author[2]{Davide Carbone}
\author[3,4]{Lamberto Rondoni}
\author[5]{Abramo Agosti}

\affil[1]{\normalsize Department of Physics, Faculty of Natural, Mathematical \& Engineering Sciences, King’s College London, London, United Kingdom}
\affil[2]{\normalsize Laboratoire de Physique de l'Ecole Normale Supérieure, 
  Université PSL, CNRS,  Sorbonne Université, Université de Paris
  Paris, France  }
\affil[3]{\normalsize Department of Mathematical Sciences, Politecnico di Torino, 10129 Torino, Italy}
\affil[4]{\normalsize INFN, Sezione di Torino, Via P. Giuria 1, 10125 Turin, Italy}
\affil[5]{\normalsize  Department of Mathematics, University of Pavia, 27100 Pavia, Italy}

\begin{document}

\title{Characterizing the Dynamics of Muscle Regeneration in the mdx Mouse Model}

\maketitle

\pagenumbering{roman}
\cleardoublepage
\cleardoublepage

\pagenumbering{arabic}
\setcounter{page}{1}

\include{Abstrac}
\include{Introduction}

\include{Chapter1}

\include{Chapter2}

\include{Conclusions}

\appendix

\include{appendix1}
\include{appendix2}

\cleardoublepage
\addcontentsline{toc}{chapter}{\bibname}

\nocite{*}

\printbibliography

\end{document}

%% file: Abstrac.tex
\section*{Abstract}
\addcontentsline{toc}{chapter}{Abstract}

In this work, we explore a model of immune response in the skeletal muscle of the \textit{mdx} mouse, the preferred animal model for studying Duchenne Muscular Dystrophy. The system studied is based on existing literature and describes the reaction of the tissue, and of the immune cells within it, to external damage, focusing specifically on the regeneration process. This work extends the analysis conducted by the authors of the original paper, highlighting interesting dynamical properties of the system, particularly in terms of long-term behavior and dependence on initial conditions. A modification is then introduced to the original equations, incorporating a diffusion term for immune cells. The dynamics of the system is numerically investigated in both the one-dimensional and two-dimensional cases. Specifically, we investigate the effect of diffusion under various boundary conditions, and the system's response to localized damage.

\keywords
{Duchenne Muscular Distrophy, Mathematical model, Dynamical system and stability, Longterm behavior, Reaction-diffusion equations}

%% file: Introduction.tex
\addcontentsline{toc}{chapter}{Introduction}


\section{Introduction: Models for Duchenne Muscular Dystrophy}

The term ``muscular dystrophy''  refers to a group of conditions, distinguishable from one another for various factors, primarily their 
underlying genetic anomaly. The main symptoms common to these disorders are progressive muscle weakness, or myopathy, and atrophy. In some cases, these manifestations are relatively mild and localised; in other forms, quality of life is severely compromised, and life expectancy is drastically reduced.
The most diffuse form of muscular dystrophy in childhood, known as Duchenne Muscular Dystrophy (DMD), affects approximately 1 in 5000 male births \cite{crisafulli}. It is a severe condition, linked to mutations in the gene encoding dystrophin on the X chromosome. For this reason, it is predominantly found in males, although milder symptoms can also occur in female carriers. The discovery of the gene responsible for DMD and the subsequent isolation of dystrophin occurred in the 1980s, allowing researchers to identify many of the molecular and genetic mechanisms underlying the disease. However, its complexity and the involvement of numerous immunological, mechanical, and chemical factors associated with the symptoms, some of which may still be unidentified, make it difficult to pinpoint specific therapeutic targets. To date, no effective treatment exists.

The predominant effects of the disease are the weakening and degeneration of skeletal muscles. Specifically, the absence of dystrophin renders the muscle fragile and unable to perform its physiological functions without incurring damage. This damage is not promptly repaired, leading to a gradual transformation into necrotic fibers and fat. Diagnosis is not always immediate, as the symptoms become clearer only when this substitution process is already in an advanced stage.
By around the age of 12, most individuals are confined to a wheelchair. As the disease advances, the onset of cardiac and respiratory complications leads to a life expectancy typically around thirty years, despite available therapies.

Beyond its primary effects on muscle tissue, the disease presents a wide variety of symptoms (pseudohypertrophy, cardiovascular disorders, hypersomnia, diabetes, retinal problems and others \cite{lapelusa}), and the anomalies are not always easy to interpret. The absence of dystrophin impacts the clinical characteristics of the disease both directly and indirectly. As a result, it is often difficult to determine whether a specific manifestation is attributable to the expression of the defective gene itself or to a cascade of events triggered by it.


Genetic mutations affecting dystrophin expression can also be found in certain animal species. Unfortunately, the clinical presentation typically differs significantly from that observed in humans. At present, dozens of different animal models exist, both naturally occurring and laboratory-generated \cite{mcGreevy}. The most evident similarities are seen in dogs; however, for practical reasons, the preferred model for study is the mouse, specifically the \textit{mdx} mouse model.
This type of mouse was first observed in the early 1980s. Its name refers to its characteristic genetic mutation and stands for ``dystrophin mutation X-linked''. Despite the X-linked inheritance and some other shared features, such as the predominant involvement of skeletal muscles, the \textit{mdx} mouse does not exhibit the profoundly debilitating symptoms seen in patients with DMD. 

Quite remarkably, the muscle histology differs substantially from that of humans: damaged fibers are repaired, and fibrosis does not occur. The damage is milder, and life expectancy is not significantly shortened. It is approximately  $25\%$ versus the $75\%$ concerning human patients \cite{mcGreevy}) . Additionally, \textit{mdx} mice are resistant to pathological intramuscular fat infiltration (IMAP) \cite{norris}.
This discrepancy can be partially attributed to fundamental differences in scale between humans and mice, which affects biomechanics, metabolism, and tissue structure. These scaling differences also influence molecular diffusion and should be taken into consideration when dealing with such mechanisms, as they may further limit the applicability of findings from mice to humans \cite{partridge}.

In healthy muscle, dystrophin is part of the dystrophin-associated protein complex (DAPC). The DAPC plays a role in signalling and serves as an essential structural element, linking the cytoskeleton, sarcolemma, and extracellular matrix (ECM). A deficiency in dystrophin results in a more fragile sarcolemma and increased cell permeability to the flow of macromolecules. The phenomenon of the ingress of foreign molecules into the cytoplasm is further exacerbated by mechanical stress.

The muscle's increased susceptibility to damage induced by contraction leads to a state of chronic inflammation, which the body cannot resolve. Muscle regeneration is compromised, and the muscle degenerates: mionecrosis, fibrosis (excessive ECM accumulation), fat infiltration, and more occur. The muscle gradually loses its elasticity, and as the disease progresses, the regenerative capacity continues to decline. Atrophy occurs after an initial phase of pseudohypertrophy, where the replacement of muscle fibers with fibrous connective tissue and fat gives the false impression of an overdeveloped muscle.

A crucial interaction in the regeneration process is that between the endogenous satellite cell population and immune system cells. The types of immune cells involved are diverse and must not only act in synergy but also with proper timing \cite{ziemkiewicz} for muscle regeneration to occur correctly. These cells display different phenotypes, both pro-inflammatory and anti-inflammatory. Both actions are necessary for an effective immune response, although it is still unclear exactly how and to what extent inflammation is beneficial to repair.

When muscle is damaged, the body reacts by triggering a series of events that lead to the infiltration of immune cells and the release of pro-inflammatory cytokines. These attract neutrophils, which, in turn, secrete cytokines and other substances, initiating an inflammatory response through phagocytosis. The presence of neutrophils also leads to an increase in macrophage concentration, which is critically important in the regeneration process. Macrophages derive from the differentiation of monocytes and perform various functions, from enzyme secretion to phagocytosis and recruitment of other immune cells. Depending on their state, macrophages display different phenotypes and functionalities. What triggers the transition from one state to another is not yet fully understood, as \textit{in vivo} differences are not rigid, and there is considerable heterogeneity.

Other key players in the immune response to damage are lymphocytes, specifically CD4$^+$ and CD8$^+$ T cells, which increase in concentration in the presence of necrosis. These cells have proven to be a good target for therapies against dystrophy, as the chronic damage to the affected muscle prolongs their presence in the tissue, further exacerbating the inflammation \cite{ziemkiewicz}. In this regard, it is important to highlight that the prolongation of the damage induces persistent inflammation, disrupting the repair mechanisms. The lymphocytes of the first family are called \textit{helper}. The lymphocytes of the second family are called \textit{cytotoxic}. The action of lymphocytes in damaged muscle is highly varied, even within the same family. It seems that lymphocytes can drive the tissue into an inflammatory or anti-inflammatory state through the secretion of various substances, influence the polarization of macrophages, and affect the proliferation of satellite cells.


Such a complex phenomenology makes mathematical modeling extraordinarily challenging, although
most desired. Mathematical models, intended to complement experimental biological studies, can in fact explore rapidly and at virtually no cost a wide variety of scenarios.
Then, on the basis of the available experimental data, one can illustrate and predict the effects of the possible mechanisms underlying numerous aspects of interest in DMD, such as muscle repair.

In some models, the focus is on post-injury muscle regeneration in healthy tissues (\cite{stephenson}, \cite{kojouharov}).
The phenomenon of muscle regeneration can also be addressed through modelling when the damaged tissue is affected by a pathological condition. The existence of an animal model with a pathophysiology in some respects similar to that of humans makes DMD amenable to this type of approach. Thanks to the availability of experimental data on \textit{mdx} mice, mathematical models can be validated.

In the literature, there are different models that focus on specific aspects of the pathological process, given the complexity of the disease, with its different cellular, molecular, and tissue dynamics.
Notable works include Ref.\ \cite{virgilio} by Virgilio, Martin, Peirce, and Blemker, as well as Ref.\ \cite{houston}, by Houston and Gutierrez \cite{houston}. The former developed an agent-based model focused on a 50-fiber muscle section in the lower limb of a mouse, analyzing regeneration following acute injury through the interaction of spatial agents and a wide range of inflammatory cells. In the latter, the focus shifted to the role of the immune system, specifically the plasticity of macrophages in the regeneration process. Various functional forms were tested to simulate different types of damage, which allows the analysis of scenarios in which damage continually fuels inflammation, providing a more realistic representation of the pathological dynamics.

In a $2009$ study \cite{dellAcqua}, Dell'Acqua and Castiglione focused on the effect of the immune system on the muscle regeneration process. The study was further developed by Jarrah, Castiglione, Evans, Grange, and Laubenbacher, resulting in additional advancements \cite{jarrah}, on which our
present work is based. 

{Our paper is organized as follows. Section 2 is devoted to a detailed analysis of the ODE model of Ref.\cite{jarrah}, 
also comparing our trajectories with those reported in that paper. In particular, we numerically investigate the long-time behavior, 
including fixed-points, 
linear stability, and the emergence of periodic regimes. 
Section 3 examines the divergence of the (asymptotically autonomous) vector field to quantify dissipativity and locate regions of phase-space contraction. Section 4 extends the model to reaction–diffusion in one spatial dimension, exploring boundary conditions and a crowding-induced reduction of effective reaction rates, and Section 5 presents the analogous two-dimensional study. We conclude by summarising implications and outlining directions for future investigations.
}

{Our main result is the following: in physiologically relevant ranges the dynamics is strongly dissipative and typically converges to a unique chronic-damage steady state,
yet specific initial immune-cell compositions or moderate diffusion under fixed-boundary constraints can drive sustained oscillations consistent with the cyclical damage–regeneration observed in mdx mice. Moreover, reaction terms dominate transport, so localised injuries remain largely confined unless additional directed mechanisms are introduced.}

%% file: Chapter1.tex
\setkeys{Gin}{draft}
\setkeys{Gin}{draft=false}

\section{Model description}

\label{chap:1}

The mathematical model from which we start  was introduced in Ref.\ \cite{jarrah}. It focuses on the role of the immune system in skeletal muscle regeneration, considering neither satellite cells or fibroblasts, nor space dependence.
It consists of six ordinary differential equations that describe the dynamics of muscle tissue composition and the concentration of immune cells over time. The system's state variables fall into two categories: immune cells and muscle tissue types.

Immune cells concentrations in number of cells per unit volume ($\mathrm{cells}/\mathrm{mm}^3$) include CD4$^+$ T lymphocytes (commonly referred to as helper lymphocytes), whose concentration is denoted by $H$, CD8$^+$ T lymphocytes (cytotoxic lymphocytes), represented by $C$, and macrophages, which are denoted by $M$.
Muscle tissue, on the other hand, is divided into three types: morphologically normal fibers ($N$), damaged fibers ($D$), and regenerating fibers ($R$). Regenerating fibers are in the process of repair. These muscle tissue variables are expressed as percentages of the total muscle composition. Due to the principle of mass conservation, one of these three variables is linearly dependent on the others. The system of equations is the following:
\begin{align}
    \label{eq:H}
    \frac{dH}{dt} &= d_H H_0 + k_1 D M - d_H H \\[0.1cm]
    \label{eq:C}
    \frac{dC}{dt} &= d_C C_0 + k_2 D H - d_C C \\[0.1cm]
    \label{eq:M}
    \frac{dM}{dt} &= d_M M_0 + k_3 D M - d_M M \\[0.1cm]
    \label{eq:N}
    \frac{dN}{dt} &= k_4 R - k_5 C N - \alpha(t) N  \\[0.1cm]
    \label{eq:D}
    \frac{dD}{dt} &= k_5 C N + \alpha(t) N - k_6 D M - d_D D \\[0.1cm]
    \label{eq:R}
    \frac{dR}{dt} &= k_6 D M + d_D D - k_4 R
\end{align}
where
\begin{equation}
    \alpha = \frac{h}{t \sigma \sqrt{2 \pi}} \exp\left(- \frac{(\ln{(t)} - m)^2}{2 \sigma^2}\right)
    \label{eq:alpha}
\end{equation}
represents a transient damage contribution, decaying in time; see Appendix \ref{app:A} for further discussion on damage modeling.

While this model excludes both the possibility of necrosis and the conversion of muscle tissue into fibrotic tissue or fat, which is a characteristic feature of the human form of the disease \cite{mareedu2021abnormal}, the resulting dynamics is quite complex, and includes numerous effects deserving careful analysis.

The evolution equations for the immune cells, \eqref{eq:H}, \eqref{eq:C}, and \eqref{eq:M}, represent their rate of change in a cubic volume with a side of $1$ \si{mm}. These equations include a source term proportional to the initial number of the corresponding cells, $H_0$, $C_0$, and $M_0$, as well as a term accounting for the death rates of the immune system cells. 
The remaining quantities describe a simplified version of the regeneration reactions. 
Specifically, $k_1 D M$ and $k_3 D M$ represent the recruitment of macrophages and helper lymphocytes in the presence of damage. These terms vanish when $D = 0\%$. The increase in cytotoxic lymphocytes is mediated by CD4$^+$ T cells in the presence of damage, as represented by $k_2 D H$.

$H_0$, $C_0$, and $M_0$ represent the physiological concentrations the system relaxes to in the absence of damage. It is noteworthy that their presence in the evolution equations implies that modifications of the initial conditions effectively amount to modifications of the dynamical law as well. An effect that should be investigated in detail.

Equations \eqref{eq:N}, \eqref{eq:D}, and \eqref{eq:R} describe the evolution of the composition of muscle tissue. Regeneration is triggered by damage, which is modeled by the function $\alpha(t)$. This function converts a fraction of morphologically normal fibers into damaged tissue, as described by the term $\alpha(t) N$.
As regeneration progresses, $k_4 R$ increases the number of normal fibers proportionally to the percentage of regenerating ones. Since macrophages clear debris from the area of inflammation, $k_6 D M$ represents the reduction of damaged fibers mediated by these immune cells. Other unspecified degradation mechanisms are captured by the term $d_D D$. Finally, cytotoxic lymphocytes contribute to muscle damage by acting on normal fibers, as described by $k_5 C N$.
The functional form of $\alpha(t)$ is designed to model the mechanical damage that arises from a multiplicative degradation process \cite{kolmogorov}. 

The optimal values of the thirteen parameters in the model are listed in Table \ref{tab:parameters} and have been estimated by Jarrah et al. \cite{jarrah} by fitting immune response data and muscle data from the literature.
\begin{table}[h!]
\centering
\renewcommand{\arraystretch}{1.1} 
\begin{tabular}{|>{\centering\arraybackslash}p{3cm}||>{\centering\arraybackslash}p{3cm}|>{\centering\arraybackslash}p{3cm}|}
\hline
\textbf{Parameter} & \textbf{Value} & \textbf{Unit} \\ \hline
$d_H$           & 0.83355             & \si{w}$^{-1}$           \\ \hline
$d_C$           & 1.61511            & \si{w}$^{-1}$             \\ \hline
$d_M$           & 0.781155           & \si{w}$^{-1}$     \\ \hline
$d_D$           & 1.34671           & \si{w}$^{-1}$             \\ \hline
$k_1$           & 0.0324139           & \si{pt\%}$^{-1}$ \si{w}$^{-1}$            \\ \hline
$k_2$           & 0.115375           & \si{pt\%}$^{-1}$ \si{w}$^{-1}$             \\ \hline
$k_3$           & 0.766576           & \si{pt\%}$^{-1}$ \si{w}$^{-1}$             \\ \hline
$k_4$           & 0.123848           & \si{w}$^{-1}$             \\ \hline
$k_5$           & 4.09948 $\times 10^{-3}$           & \si{cell}$^{-1}$ \si{w}$^{-1}$             \\ \hline
$k_6$           & 3.23097 $\times 10^{-4}$           & \si{cell}$^{-1}$ \si{w^{-1}}            \\ \hline
$\sigma$           & 2.92815           & -             \\ \hline
$m$           & 4.22686           & -             \\ \hline
$h$           & 0.511657           & -             \\ \hline
\end{tabular}
\caption{Parameters values as estimated in the reference paper. Time is measured in weeks.}
\label{tab:parameters}
\end{table}
Furthermore, the article suggests a standard initial configuration for the system, where in a volume of $1$ \si{mm^3} $H_0 = 0$ cells/mm$^3$, $C_0 = 4$ cells/mm$^3$, and $M_0 = 400$ cells/mm$^3$ (these initial counts align with the available data from wild-type mice) and excludes the possibility of damage during gestation, so that the muscle at birth is entirely composed of normal fibers ($N_0 = 100\%$, $D_0 = 0\%$ and $R_0 = 0\%$).
In the remainder of this discussion, we will refer to the system initialized in this way, with parameters set according to the specifications of the paper, as the ``default system''.

To compute the evolution of the system, the explicit Euler method was employed. This method offers a straightforward and simple approach to representing the system's dynamics. However, it is conditionally stable, necessitating careful selection of the time step $dt$ to ensure numerical stability.
A comparison with the evolution computed using the fourth-order {Runge-Kutta} method in SciPy shows results that closely align with those of the Euler method when $dt$ is set to $0.01$ weeks. This value was then adopted as the default in all subsequent simulations, unless otherwise specified. Nevertheless,
this constitutes a delicate issue, because using the same parameters of the cited article, the trends we observe over a period of $15$ \si{weeks} slightly differ from the corresponding plots in Ref.\ \cite{jarrah}. Since our Euler and Runge-Kutta implementations produce nearly identical results, the discrepancies with \cite{jarrah} cannot be merely attributed to numerical accuracy. In fact, slight variations of the initial conditions lead to plots that agree more closely with those in Ref.\ \cite{jarrah}. This reveals that the dynamics is quite sensible, but in non obvious fashions.

\subsection{Analysis of default system}

One aspect of interest, not addressed in Ref.\ \cite{jarrah} is the long-term behavior of the dynamics. Although asymptotic trends may appear too remote in time to be relevant within the life span of mice, they bear information on the evolutions allowed by the model, and these may also  
be reached within short times, under suitable conditions. Then, medical treatment could be aimed at creating favorable conditions.

\begin{figure}[ht!]
    \centering
    \includegraphics[width=\textwidth, trim=0 0 0 30, clip]{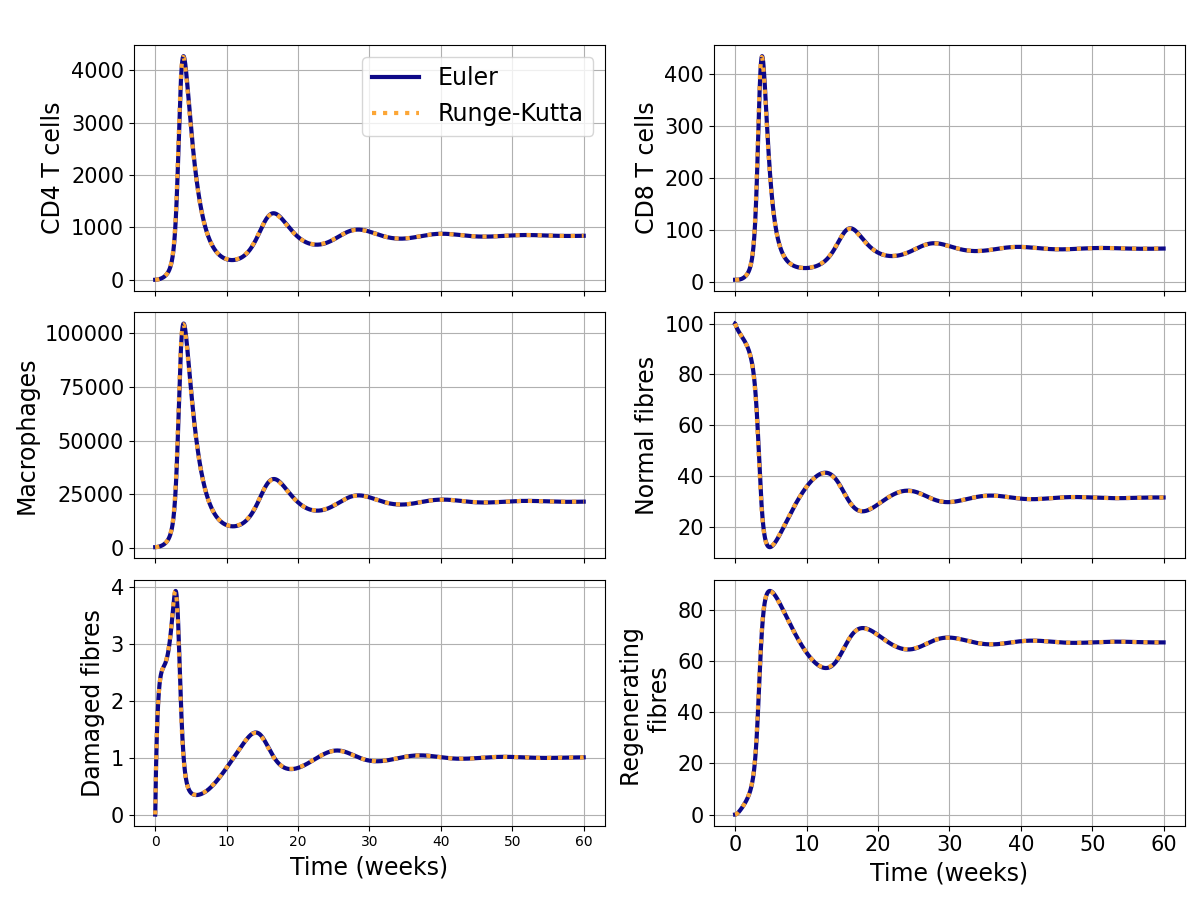}
        \caption{Evolution of immune cell concentration (\si{cells/mm}$^3$) and muscle tissue composition ($\%$) over $60$ \si{weeks}.}
        \label{fig:E_RK_60}
\end{figure}

We found that the selected initial conditions lead to a steady state which is a fixed point, with just under $70\%$ of regenerating fibers and only about $1\%$, of damaged fibers. Figure \ref{fig:E_RK_60} presents an example of the evolution over $60$ \si{weeks}. After a first peak in the concentration of immune cells, which corresponds to a massive immune response triggered by damage, several less pronounced peaks follow. Similar oscillations are observed in the percentages of normal, damaged and regenerating fibers. In the steady state, the tissue has not fully recovered its healthy state. Instead, a chronic immune response develops within it, maintaining a stable fraction of damaged fibers through a continuous regeneration process. 

A systematic analysis of the fixed points has then been carried out numerically, considering very long time spans. As physical ranges for the
initial conditions for N, D, and R we took the interval $[0 , 100]$, while for the immune cell concentrations the upper limit has been set to the value of the first peak reached in the default model.
A single physically meaningful and linearly stable fixed point was found as will be shown later.

The basin of attraction in the full phase space cannot be directly studied, because of the explicit 
dependence of the model on $H_0$, $C_0$ and $M_0$. Changes in the initial configuration imply changes in the dynamical law, hence in the stationary state. 

\begin{figure}[ht!]
    \centering
    \begin{minipage}[b]{0.48\textwidth}
        \centering
        \includegraphics[width=\textwidth, trim=140 55 95 60, clip]{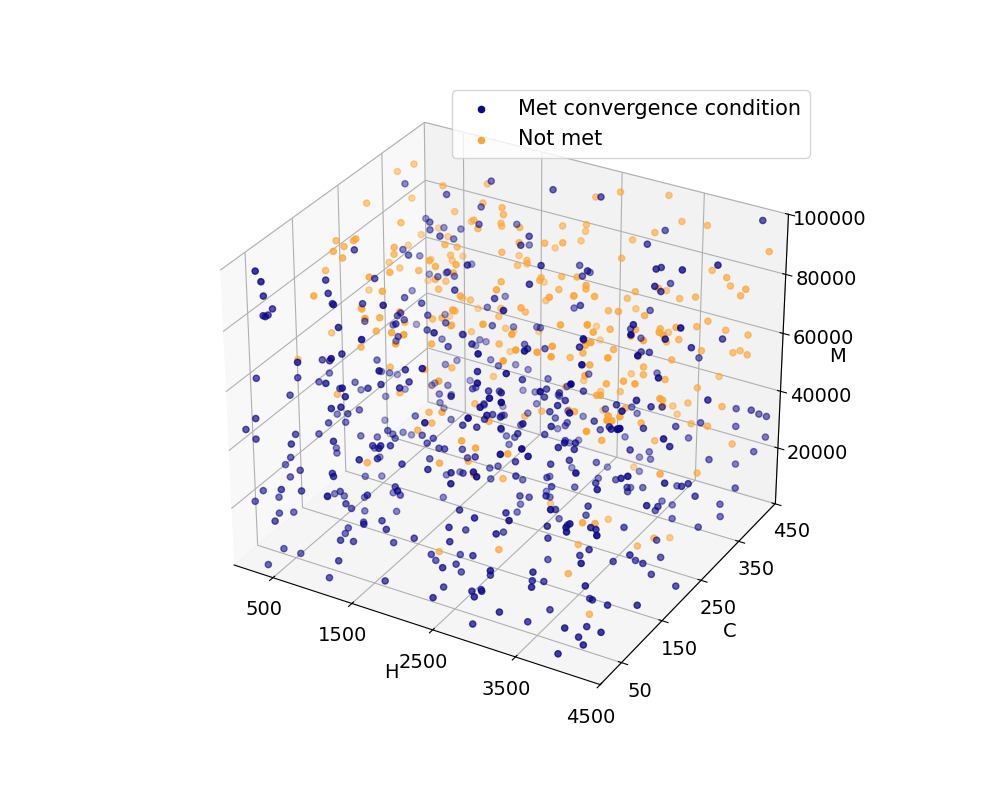}
        \captionsetup{labelsep=none}
    \end{minipage}
    \hfill
    \begin{minipage}[b]{0.48\textwidth}
        \centering
        \includegraphics[width=\textwidth, trim=140 55 95 60, clip]{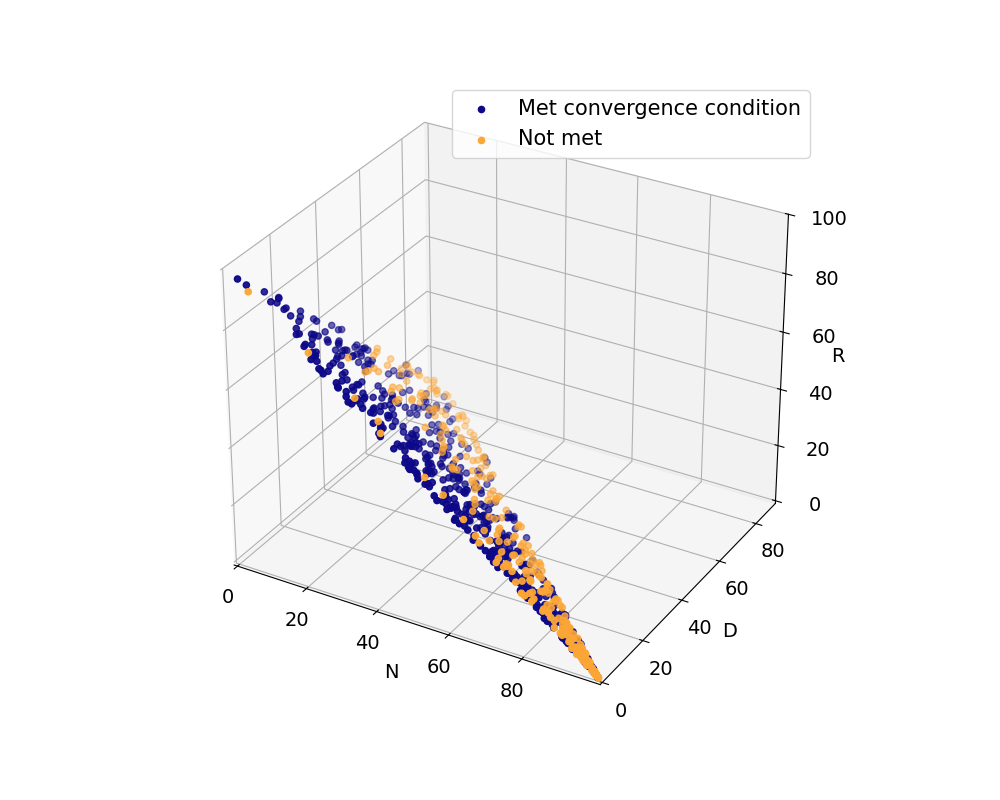} 
        \captionsetup{labelsep=none}
    \end{minipage}
    \caption{Visualisation in the $HCM$ (left) and $NDR$ (right) reduced spaces. The two colors distinguish the samples of initial configurations: orange represents cases where the system satisfies the convergence condition within $300$ \si{weeks}, while blue represents those for which it does not.}
    \label{fig:HCM_vs_NDR}
\end{figure}

In Ref.\ \cite{jarrah}, this aspect is only partially investigated and it is concluded that the system's qualitative behavior is not sensitive to the initial conditions. However, we found that this is the case as long as the initial conditions do not differ significantly from the default ones.
%
More generally, exploring the full five-dimensional phase space and searching for the asymptotic states is rather costly. However, one can gather useful information by running a relatively large number of initial conditions 
over reasonably long times, and testing a stationarity condition for each of the corresponding trajectories.
We expect a behavior qualitatively similar to that of the default system, hence a
recognizable convergence trend well before $100$ weeks. We then visualize the initial conditions in the reduced space of three out of the six variables, and divide them into two categories: those verifying the convergence condition at a certain time, and those that do not verify it.
In particular, we checked whether the configuration of the system at $271$, $283$ and $300$ weeks remained unchanged or not, within a certain tolerance, cf.\ Figure \ref{fig:HCM_vs_NDR}) for the triples $(HCM)$ and $(NDR)$.

For some of the initial configurations that failed to meet the convergence condition, the evolution differs substantially from that of Ref.\ \cite{jarrah}. Either convergence to a fixed point is significantly delayed, or a completely different asymptotic regime is present. To systematically investigate such an issue, we considered the default condition, modifying one coordinate at a time, considering that modifying one variable of the muscle composition, the others must also be modified to ensure that mass is conserved. This way we have isolated and highlighted the influence of each individual coordinate, and thus each state variable.

Before doing that, let us consider the autonomous set of ODEs, because their time dependent term vanishes in the long time limit.
%
The process of muscle regeneration is triggered by damage. In our setup, we exclude the presence of damaged fibers at birth, {\em i.e.}\ $D_0 = D(t_0)=0\%$. 
However, the term $k_5CN$ in Equations \eqref{eq:N} and \eqref{eq:D}, combined with a non-zero initial concentration of CD8$^+$ T cells, allows an evolution of the variables $D$ and $N$ even in the absence of $\alpha(t)$. Normal fibers are partially converted into damaged fibers in the presence of cytotoxic lymphocytes, thereby activating the dynamics, 
although there is no guarantee that  autonomous and non-autonomous dynamics  will share the same fixed point. Even if $\alpha(t)$ decays over time, during the transient phase it could steer the system towards a different basin of attraction before vanishing, leading to a shift that is not observed in the autonomous system. 
Furthermore, it is not possible to systematically directly compare the effect of $\alpha(t)N$ with that of the other terms governing the evolution of the percentage of normal and damaged fibers, because the value of these terms  depends on the initial conditions.

Nevertheless,
a test performed on the default system showed a decay of $\alpha(t)$, making it rapidly negligible relative to the other terms. 
Since the damage term appears to have an effect limited to the first few weeks of evolution, it is unlikely that different initial conditions would make it strong enough to alter the steady state. Therefore, in all subsequent simulations, it is assumed \textit{a priori} that the autonomous and non-autonomous systems share the same final equilibrium configuration regardless of the initialization. Some tests were conducted to validate this assumption. The search for the convergence time, which we will call $t_c$, can therefore make use of this additional condition.


Indeed, the assumption that the autonomous and non-autonomous systems become indistinguishable at steady state provides a more immediate criterion for testing convergence: when the trends with and without $\alpha(t)$ become sufficiently close (within a given tolerance) and remain so for the rest of the evolution, the configuration is likely to be similar to that of the steady state. 
We then implemented a convergence criterion to halt the evolution:
the difference between the value of each state variable at time $t$ and $t+1$ must be smaller than a given tolerance.
Letting $a$ be the value of any state variable at time $t$ and $b$ its value at the next time step, the following inequality must hold
\begin{equation}
    \frac{|a - b|}{\max(a, b)} \leq tol,
    \label{eq:rel_tol}
\end{equation}
where $tol$ defines an arbitrary tolerance. Equation \eqref{eq:rel_tol} is also used to assess the similarity between the autonomous and non-autonomous configurations.

This way, we found that convergence is  problematic to assess when the initial concentration of macrophages is varied; if inequality \eqref{eq:rel_tol} is never satisfied after $400$ weeks, we record the final value, being that time horizon already wll beyond the average life of the mouse mdx. We did not encounter a similar issue for the other state variables.
For the root-finding problem, we used standard routines of \texttt{scikit-learn} package. Then,
due to the numerical limitations of the solvers, with different parameters and setups, for identifying fixed points in the default system, 
we relied on  
the convergence times procedure to examine how the system's configuration evolves, under the assumption that our test actually identifies a steady state.
%
Then, one of the state variables is selected and starting from the initial configuration $v_i$, 
the system's evolution is computed over a span of $400$ weeks, with states recorded at both $300$ and $400$ weeks. The values of each state variable at these checkpoints are then displayed as a function of the initial condition $v_i$.

As expected, when the selected variable pertains to the composition of muscle tissue, the steady-state remains unchanged.
A different scenario arises for variables associated with immune cell concentrations. Since $H_0$ and $C_0$ appear explicitly in the equations of the model, the fixed point shifts position in phase space when these conditions are altered. The variation appears continuous, exhibiting a smooth transition as the aforementioned initial conditions change. The configurations at $300$ and $400$ weeks overlap, which confirms that the displayed configurations represent a fixed point of the system. The linear stability analysis performed around the system state reached at very long times, treated as a fixed point under the assumption that the dynamics have already converged to a steady regime, further confirms that these equilibria are stable for all the initial conditions tested.

\begin{figure}[ht!]
    \centering
    \includegraphics[width=1.0\linewidth, trim=0 0 0 10, clip]{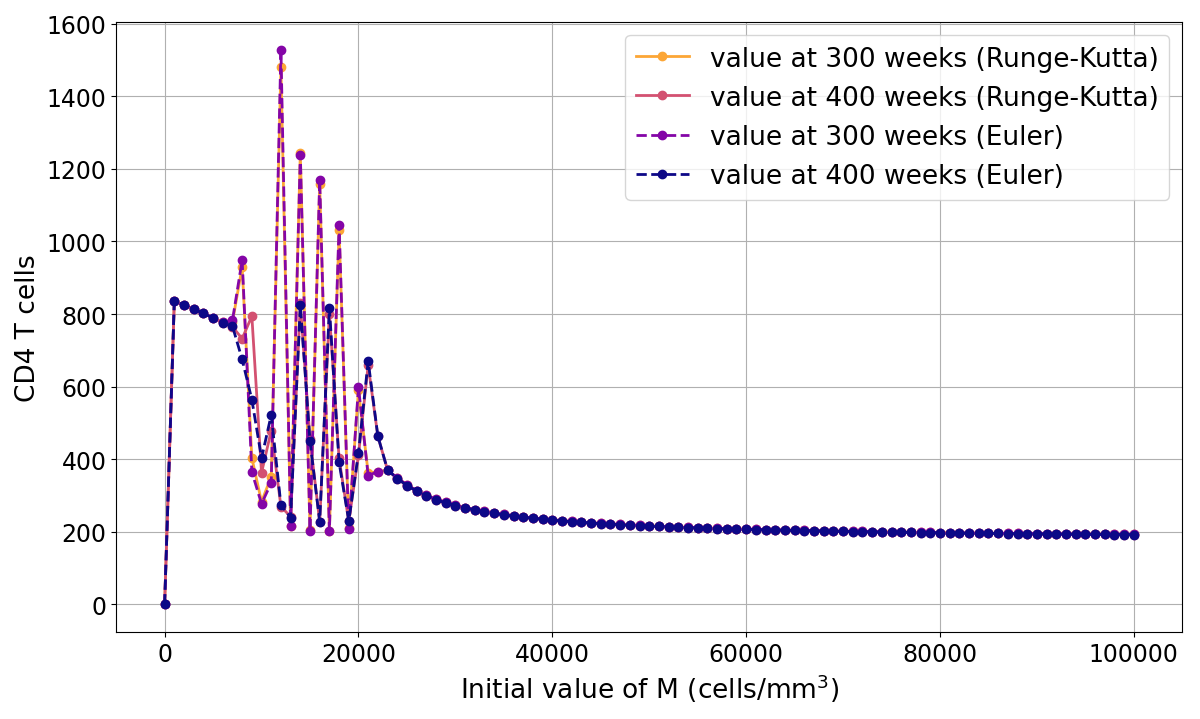}
        \caption{Behavior of the system at $300$ and $400$ weeks as the initial concentration of macrophages varies within the range $[0, 100000]$ \si{cells/mm}$^3$.}
        \label{fig:fp_M}
\end{figure}

The behavior of the system state at $300$ and $400$ weeks, as the initial condition for macrophages varies, appears to be rather complex, with  Euler and Runge-Kutta methods yielding the same. In the first place, at small values of $M_0$, 
the state variables alleged
steady state values appear discontinuous as functions of the initial condition, cf.\ Figure \ref{fig:fp_M}. For a wide range of initial concentration values, we find a fixed point whose position in phase space varies continuously with changes in $M_0$. For very low concentrations a kind of discontinuity arises.
This depends on the choice of the initial condition for the helper lymphocytes.
The combination of $M_0 = 0$ \si{cells/mm}$^3$ and $H_0 = 0$ \si{cells/mm}$^3$ prevents any immune cell concentrations from evolving, as their update at the first time step remains zero, completely halting their dynamics. This explains the apparent discontinuity
illustrated in the plots.

\begin{figure}[ht!]
    \centering
\includegraphics[width=1.0\linewidth]{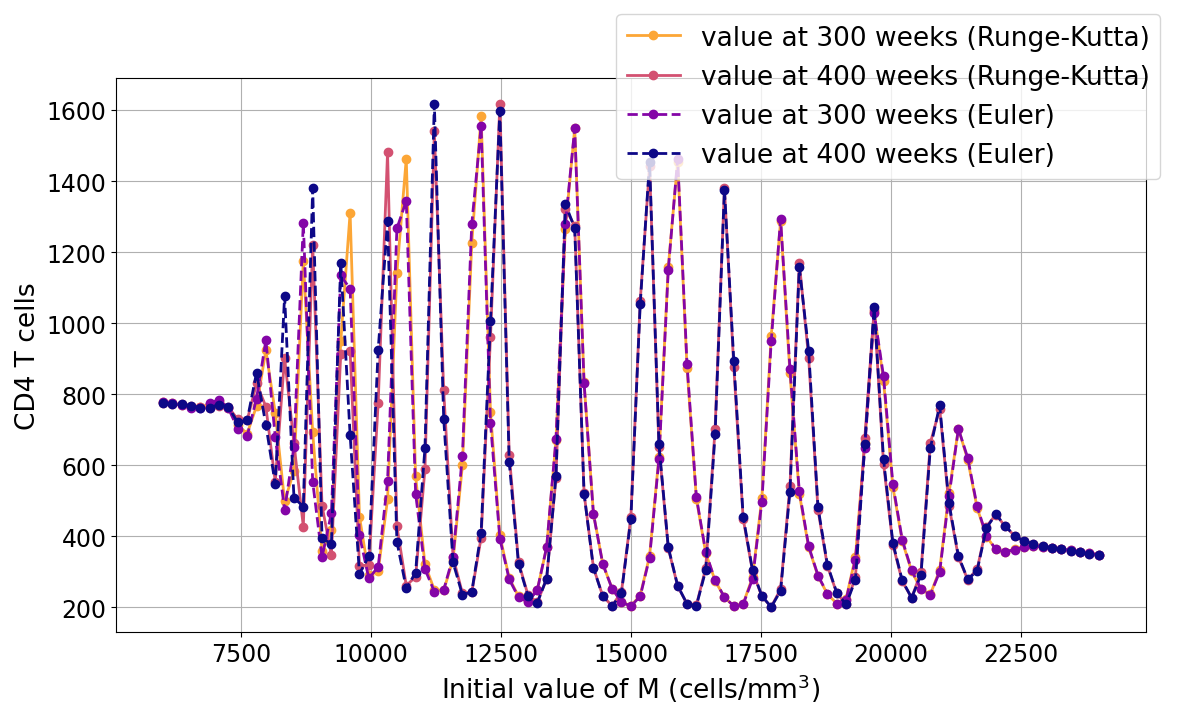}
        \caption{behavior of the system at $300$ and $400$ weeks with equally spaced $M_0$ in $[6000, 24000]$ \si{cells/mm}$^3$. The values at $300$ and $400$ weeks are different. However, the results obtained with Euler and Runge-Kutta are
        quite close. A fixed point has not been approached.}
        \label{fig:issue_M}
\end{figure}

Furthermore, consider Figure \ref{fig:issue_M}, that provides a detailed view of the range $[6000, 24000]$ \si{cells/mm}$^3$ from the first plot in Figure \ref{fig:fp_M}. At least for many initial conditions within this range, the Runge-Kutta and Euler methods produce similar results, but at $300$ and $400$ weeks significant differences appear. At the same time, the accurate evaluation of the convergence time as $M_0$ varies becomes harder and harder. 
There are two possible explanations: either the system will eventually converge, but requiring exceedingly long times, or the system does not converge at all to a fixed point, allowing for different evolutions under minimally different conditions.
Assuming that a fixed point has been correctly identified in the other regions of the phase space,  
one may attempt a
linear stability analysis, to gain insight in the possible asymptotic regimes. One finds that small $M_0$ leads to all negative eigenvalues for the Jacobian matrix, which then implies that a possible fixed point is stable.
As $M_0$ increases, around $7000$ \si{cells/mm}$^3$, some eigenvalues attain zero real parts. Further increasing $M_0$ to approximately $22000$ \si{cells/mm}$^3$ restores negative values for all the eigenvalues.

\begin{figure}[ht!]
    \centering
    \includegraphics[width=1.0\linewidth, trim=0 0 1
    0 0, clip]{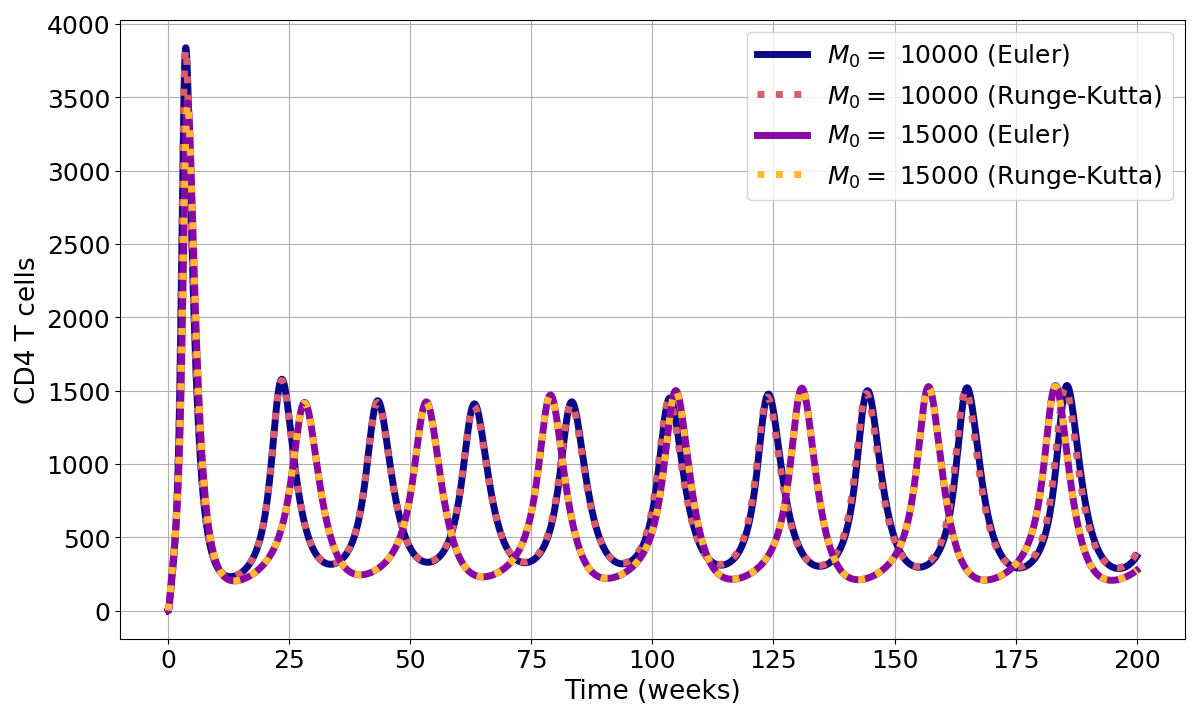}
        \caption{Evolution of the system over $200$ weeks starting from initial conditions $\{ 0,4, 10000, 100, 0, 0\}$ and $\{0, 4, 15000, 100, 0, 0\}$. }
        \label{fig:periodic_M}
\end{figure}

To understand the evolution of the model in this particular scenario, we select two different initial conditions within the range $[7000,22000]$ \si{cells/mm}$^3$ and integrate the dynamics over a long time. The results, obtained using $M_0 = 10000$ \si{cells/mm}$^3$ and $M_0 = 15000$ \si{cells/mm}$^3$, are presented in Figure \ref{fig:periodic_M},
and reveal periodic behavior. This fact is significant because to the best of our knowledge it had not been observed before, and because it aligns with  
experimental observations describing cycles of damage and regeneration in the \textit{mdx} mouse \cite{Yokota}. The trends indicate that the initial macrophage concentration influences the cycle period, and we found that this period increases monotonically with growing $M_0$, in the range from 8000 \si{cells/mm}$^3$ to $20000$ \si{cells/mm}$^3$.

\section{Vector field divergence}


An interesting indicator of the dynamics of a dynamical system is the divergence of its vector field.
Since the system of Ref. \cite{jarrah} 
is asymptotically autonomous, it may be useful to consider the model without the damage term to assess the evolution of phase space volumes, once the time dependent term has decayed, something that occurs in relatively short times.
In this case, it is also conveninet to remove Equation \eqref{eq:R}, given the conservation of mass that yields $R = 100-N-D$.
Let us denote the associated vector field as $\mathcal{N}$, in $\mathbb{R}^5$. Its divergence is
\begin{equation}
 \Lambda = \nabla \cdot \mathcal{N} = - d_H - d_C - d_M + k_3 D - d_D - k_4 - k_5 C - k_6 M.
    \label{eq:non_aut_div}
\end{equation}
The state variables represent non-negative physical quantities (cf. Appendix \ref{app:pos} for a simple proof), which means that in the physically relevant region of phase space most terms in the divergence expansion are positive semidefinite. In fact, the only term in Equation \eqref{eq:non_aut_div} that can lead to an expansion of volumes is  $k_3 D$.
For this to occur, the percentage of damaged fibers must be quite high.
The minimum percentage of damaged fibers $\Tilde{D}$ at which the divergence changes sign is the 
following \footnote{In the divergence formula, all the terms are negative except $k_3 D$. Therefore, to make the divergence positive, $D$ must be sufficiently large to offset the (negative) sum of all the other contributions. Some of these contributions are fixed because they depend only on the estimated parameters, whereas $k_5 C$ and $k_6 M$ can become very large if the concentrations of $C$ and $M$ are high. The limiting case is $C=0$ and $M=0$: under these conditions, one can enforce the divergence to be zero and obtain this lower bound for $D$.}
\begin{equation}
    \widetilde{D}  = \frac{1}{k_3} \big( d_H + d_C + d_M + d_D + k_4 \big) 
     \approx 6.13 
    \label{eq:D_tilde}
\end{equation}     
This value represents a lower bound below which the physical system is strictly dissipative. In the evolution of the default system, such a high peak of damaged muscle is never actually reached. In general, for the other tested initial conditions, the steady-state value of damaged fibers remains well below this threshold, confirming that the long-term trajectory of the system is confined to regions of dimension lower than that of the phase space. 

Since the divergence $\Lambda$ depends only on three state variables, we can fix the value of one of them and construct two-dimensional heat maps that highlight the dependence on the other two.
To select the value ranges for displaying the heat maps, we can use the results obtained from the steady-state analysis while varying the initial conditions one variable at a time. Considering asymptotic behaviors ($\alpha(t) \approx 0$), we focus on the steady state values of the variables $C$, $M$, and $D$, from which $\Lambda$ depends.
It is worth noting that varying only one initial condition at a time limits the range of accessible steady states. Therefore, for other initial configurations, the final concentrations of cytotoxic lymphocytes and macrophages may exceed the limits we found with our simulations, and the percentage of damaged fibers could also be higher at steady state.

\begin{figure}[ht!]
    \centering
    \includegraphics[width=1.0\linewidth, trim=0 0 40 40, clip]{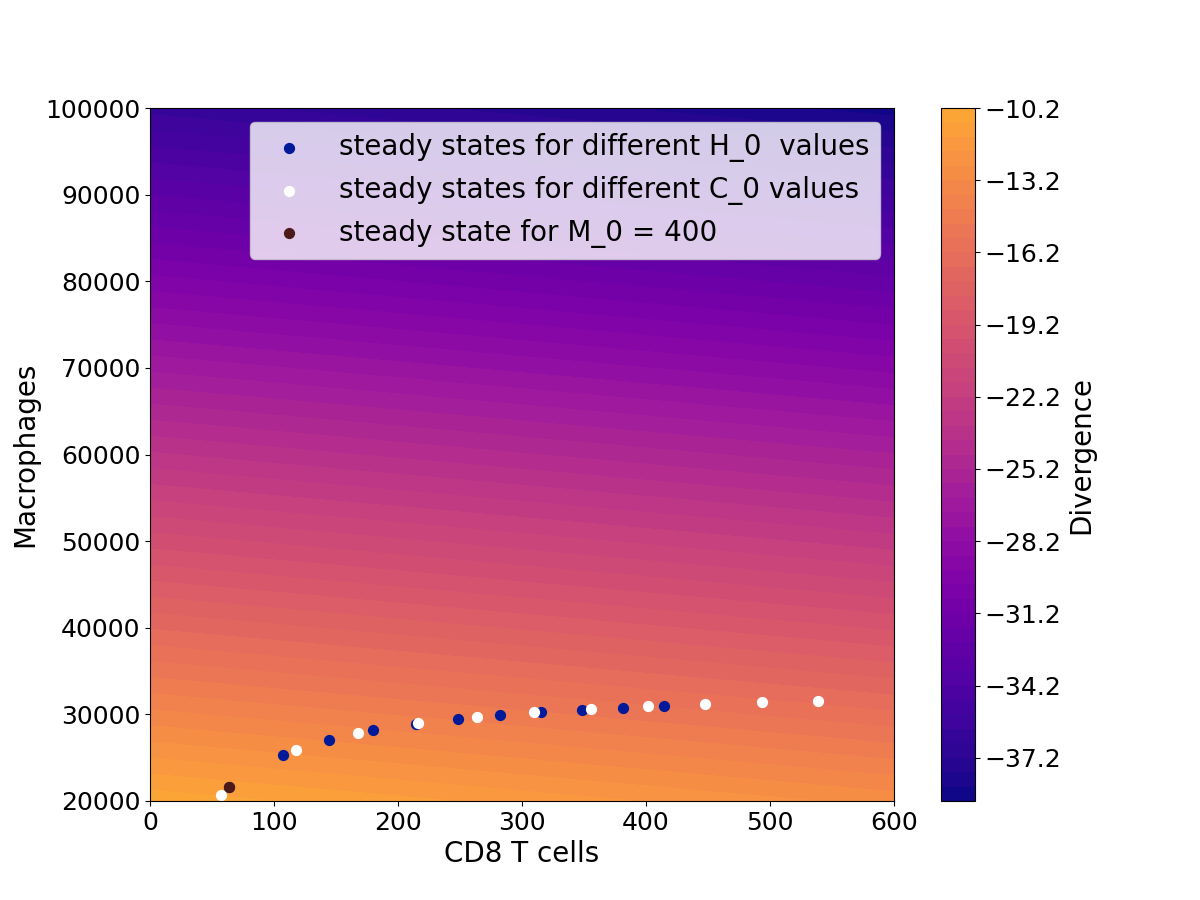}
        \caption{Heat map of the divergence of the vector field for fixed $D=1.00\%$ and different concentration values for CD8$^+$ T lymphocytes and macrophages (\si{cells/mm}$^3$). In the highlighted region, the divergence is strongly negative. The positions of some steady states are indicated by colored dots. The blue points correspond to a sample of fixed points obtained by varying only the initial condition of the helper lymphocytes, the white points correspond to a variation in the initial coordinate of the cytotoxic lymphocytes, and finally, the brown point represents the steady state of the system initialized as $\{0,4,400,100,0,0\}$.}
        \label{fig:div_D_1}
\end{figure}

Here, we only show the heat map of the divergence evaluated for a fixed percentage of damaged fibers (Figure \ref{fig:div_D_1}). This choice is linked to the fact that the value of $D$ at steady state appears to remain fairly stable despite variations in the initial conditions. In fact, for many of the tests carried out, the steady-state percentage of damaged fibers is approximately $1.00\%$, with variations affecting the second or third decimal place. This makes it possible to visualise different steady-state configurations found in the previous analysis on the heat map because, of the three variables on which divergence depends, the one related to damaged fibers can be considered fixed.

Such steady states originate from the analysis of the system as the initial condition on helper lymphocytes, cytotoxic lymphocytes, or macrophages varies, and they are color-coded accordingly. In the first two cases, virtually all equilibrium configurations satisfy the condition on the steady-state value of damaged fibers, whereas for variable $M_0$, the percentage of damaged fibers is less stable. As a result, a single point for $M_0 = 400$ \si{cells/mm}$^3$ is displayed on the map.

In this region of the reduced space, the system is everywhere dissipative. The steady-state values shown are obtained from an initial condition that differs from the default one in just one coordinate. However, the coordinate is not always the same. Despite this, the system appears to be converging towards the same manifold in phase space.

When the variables kept fixed are the concentrations $C$ and $M$, we do not observe any preferred values for which many steady states are found, as was the case in the scenario just analized. In all the scenarios considered, whose results are omitted here, the divergence of the vector field remains negative. From this, we conclude that all the previously identified steady states exist in regions where phase space volumes contract and must be even non-linearly stable.



Although the steady state appears to be located in a dissipative region of phase space for all the initial conditions considered, nothing guarantees that the system does not explore zones where volumes locally expand during the transient phase of the evolution. When studying the steady state as the initial conditions vary, it was not necessary to analize the detailed evolution of the system. It may be interesting to observe how the divergence of the vector field changes over time, starting from a limited number of specific initializations.

As previously observed, in the default system the percentage of damaged fibers never exceeds the threshold $\Tilde{D}$ calculated earlier \eqref{eq:D_tilde}, and indeed, we find that the divergence remains consistently negative over time.
This is not true for all explored initial conditions. For example, for certain values of the initial concentration of cytotoxic lymphocytes (e.g. $C_0 = 150$ \si{cells/mm}$^3$), the divergence initially takes on positive values, but then stabilizes at a steady state with negative divergence.

\begin{figure}[ht!]
    \centering
    \includegraphics[width=1.0\linewidth]{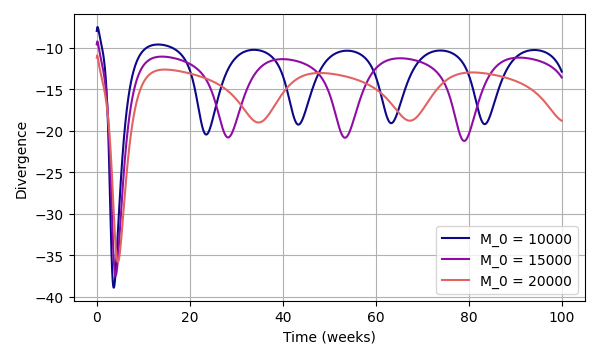}
\caption{Evolution of the vector field's divergence, starting from different initial conditions for the variable $M$. These conditions belong to the range of values for which a periodic long-term evolution trend has been observed. Specifically, the chosen values are $M_0 = 10000$ \si{cells/mm}$^3$, $M_0 = 15000$ \si{cells/mm}$^3$, and $M_0 = 20000$ \si{cells/mm}$^3$.}
        \label{fig:div_Mi}
\end{figure}
Another situation of interest is when the system stops converging to a stable fixed point and the long-term behavior becomes periodic. Figure \ref{fig:div_Mi} shows the evolution of the divergence for three different initial conditions on the $M$ coordinate, for which the system exhibits periodic behavior. In these scenarios, $\Lambda$ is strongly negative not only at the steady state but throughout the evolution, and in particular it never becomes positive, either in the long term or during the transient phase. This practically excludes complex dynamics in time, although not in the dependence on the model parameters.

In general, the divergence of the vector field can be negative in some areas of phase space and positive in some others, causing volumes to contract and expand depending on where they are located. However, from the analysis carried out, it appears that most of the physically accessible phase space tends to be characterized by quite strong volume contraction, since it proceeds at the exponential rate 
$\Lambda$ that is sensibly negative.

\noindent

%% file: Chapter2.tex
\section{Reaction and diffusion}
\label{chap:2}

\setkeys{Gin}{draft}
\setkeys{Gin}{draft=false}

The model analyzed so far describes the evolution of immune cell concentrations and muscle tissue composition through a system of ordinary differential equations, which takes into account the interactions between the variables, assuming that these quantities are uniformly distributed in space. This condition is clearly not met in skeletal muscle tissue.
In fact, immune cells move both via active migration and diffusion, spreading from areas of high to low concentration. For a more realistic description of the process, we modified the model of Ref.\ \cite{jarrah} to incorporate a term that accounts for immune cell movement during muscle damage and regeneration.


\subsection{One space dimension}
Introducing explicit space dependence can be done in different fashions. If space is continuous, one obtains a system of PDEs, but one may as well consider discrete space, which is in any case needed also to solve numerically the PDEs, and then one continues to deal with ODEs. 
This has various advantages. Mathematically, existence and uniqueness of solutions of equations with smooth vector fields in  compact spaces are guaranteed. Moreover, the cells of our interest are relatively large, and subdividing space in  domains sufficiently large to contain a large number of them seems to be a suitable approach. As a matter of fact there are numerous successfull discrete models of spatially extended systems, including cellular automata, and lattices with which numerous kinds of transport can be described, cf. \cite{chopard1998modeling,wolfram1983statistical,chopard2002cellular}.

Therefore, in one dimension, we consider a row of $N_c \ge 3$ cubes with side length of $\delta x = 1$ \si{mm}, and we use the explicit Euler integration scheme to account also for the diffusion between neighboring cubes. The behavior of the boundary cubes is then determined by the boundary conditions. 
Three different types of boundaries have been implemented: periodic boundary conditions, edge cells with a fixed concentration of immune cells,
\footnote{This condition simulates the presence of immune cell reservoirs at the boundaries.} 
and evolving edge cells, with no immune cells flux from the outside \footnote{This type of condition can be interpreted as representing a physical barrier, such as a membrane or connective tissue.}, see Appendix \ref{app:B} for further details.
%
%
%
The model then aims to describe
the action of three different types of immune cells: helper lymphocytes, cytotoxic lymphocytes, and macrophages. These cells do not serve the same functions and are characterized by different sizes and shapes. This clearly affects the way they can diffuse. Despite their clear
differences, the analysis of the system with diffusion has been simplified
by assuming a single diffusion coefficient, $\mathscr{D}$, common to all
three cell types ($D_H = D_C = D_M$). Indeed, analogous diffusion
properties may be enjoyed also by different cells, for different reasons.
Moreover, to the best of our knowledge, systematic experimental
measurements of immune-cell diffusion within skeletal muscle are lacking,
which prevents a more precise, cell-specific characterization; in the
absence of such data, there is no compelling reason to adopt different
coefficients for the three populations. The values assigned to
$\mathscr{D}$ in our simulations may also differ significantly from the
actual coefficients that characterize diffusion within muscle tissue.
Therefore, we explored diffusion coefficients across three different orders of magnitude: $10^0 = 1$ \si{mm}$^2$/\si{week}, $10^{-1}$ \si{mm}$^2$/\si{week}, and $10^{-2}$ \si{mm}$^2$/\si{week}. The only estimate found in the literature \cite{zlobina} refers solely to macrophage diffusion. It suggests lower plausible values, 
but this is a single result, obtained from \textit{in vitro} data. Adding that the real dynamics is not one-dimensional, and also that has an effect on diffusion, it remains worth considering different values of $\mathscr{D}$.

In all our simulations, $dx = 1$ \si{mm}, and the parameter used to study variations in system behavior is $\mathscr{D}$. However, considering a smaller cubic cell size while keeping the diffusion coefficient fixed has the same effect as modifying the diffusion coefficient in a system where $dx$ remains unchanged. In other words, by choosing a finer discretisation, we automatically obtain the same system evolutions, but with lower diffusion coefficients that are compatible with the few available data.



We begin with the simplest possible situation: homogeneous initial distribution of cells. Under these circumstances, neither periodic boundary conditions nor those where the boundary cells evolve like the internal ones but with zero immune cell influx from the outside are of interest. The only interesting scenario is with fixed boundary conditions.


Initially, we consider a row comprising a sufficiently large number of cells and assess the effect of diffusion for different values of the coefficient $\mathscr{D}$. Naturally, the effect of the fixed-concentration boundary will gradually spread from the outermost cells to the more central ones in the row, but the system will maintain symmetry around the center, if the two boundaries are kept at the same values.
Let us begin
\begin{figure}[ht!]
    \centering
\includegraphics[width=1.0\linewidth]{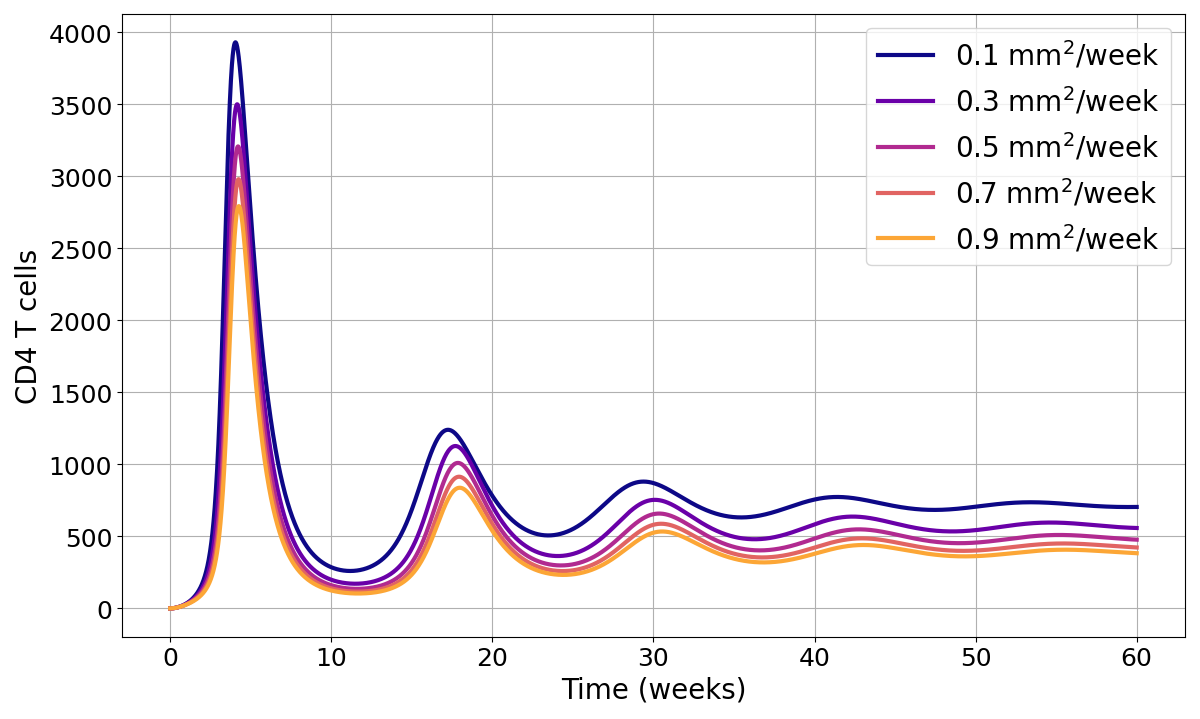}
        \caption{Evolution of the concentration of helper lymphocytes (\si{cells/mm}$^3$) in the second cell of a row of $50$ cells. Different values of the diffusion coefficient are tested.}
        \label{fig:default_diff}
\end{figure}
with the second cell of the row, which represents the site where the effect of the fixed concentration boundary is first and most strongly perceived.
 
In Figure \ref{fig:default_diff}, we can observe the evolution over a period of $60$ weeks as the diffusion coefficient $\mathscr{D}$ varies. The presence of diffusion alters the system's steady state, but does not change the qualitative trend of the state variables over time. The peaks in immune cell concentration are lower, but the convergence times remain largely unchanged.

\begin{figure}[ht!]
    \centering
    \includegraphics[width=1.0\linewidth]{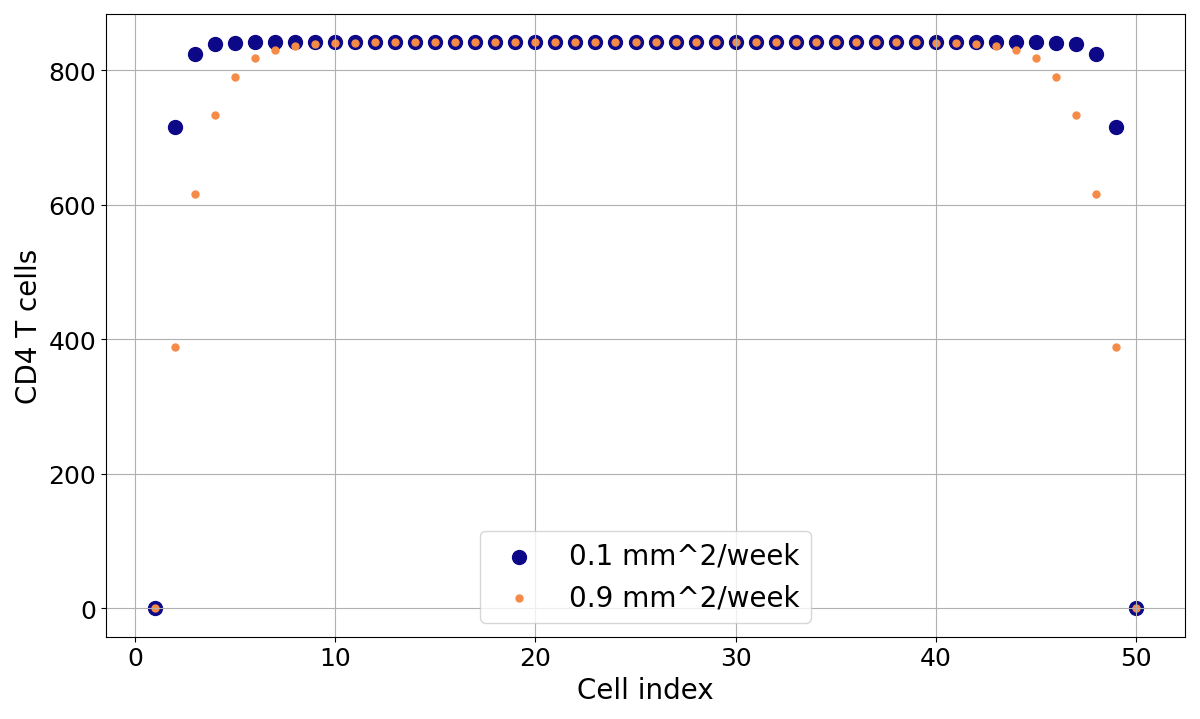}
        \caption{Helper lymphocyte concentration per cubic millimetre at $150$ weeks in all cells of a row of $50$. The results correspond to two different simulations, each characterized by a different diffusion coefficient.}
        \label{fig:H_steady_diff}
\end{figure}

A similar evolution can also be observed in the other cells of the row, although the edge effect is clearly weaker in the central positions. To assess the extent of this dependence, the configuration in all the cells can be visualized at very long times, as shown in Figure \ref{fig:H_steady_diff}.
From this comparison, we observe that in fact the diffusion effect remains limited to a small number (depending on the diffusion coefficient) of cells near the edge. It is evident that, for a significant portion of the central cells, the steady-state configuration is the same as without diffusion. 
This indicates that, in our model, the reaction process dominates over diffusion; beyond a certain distance from the boundary, cells remain unaffected by the boundary conditions.

Given the result obtained with rows composed of a high number of cells, we can reduce $N_c$ in the hope of making the effects of the fixed-concentration boundary more evident.
For instance, let us consider a row composed of $7$ cells, all once again initialized with the same conditions suggested in the paper.

\begin{figure}[H]
    \centering
    \includegraphics[width=1.0\linewidth]{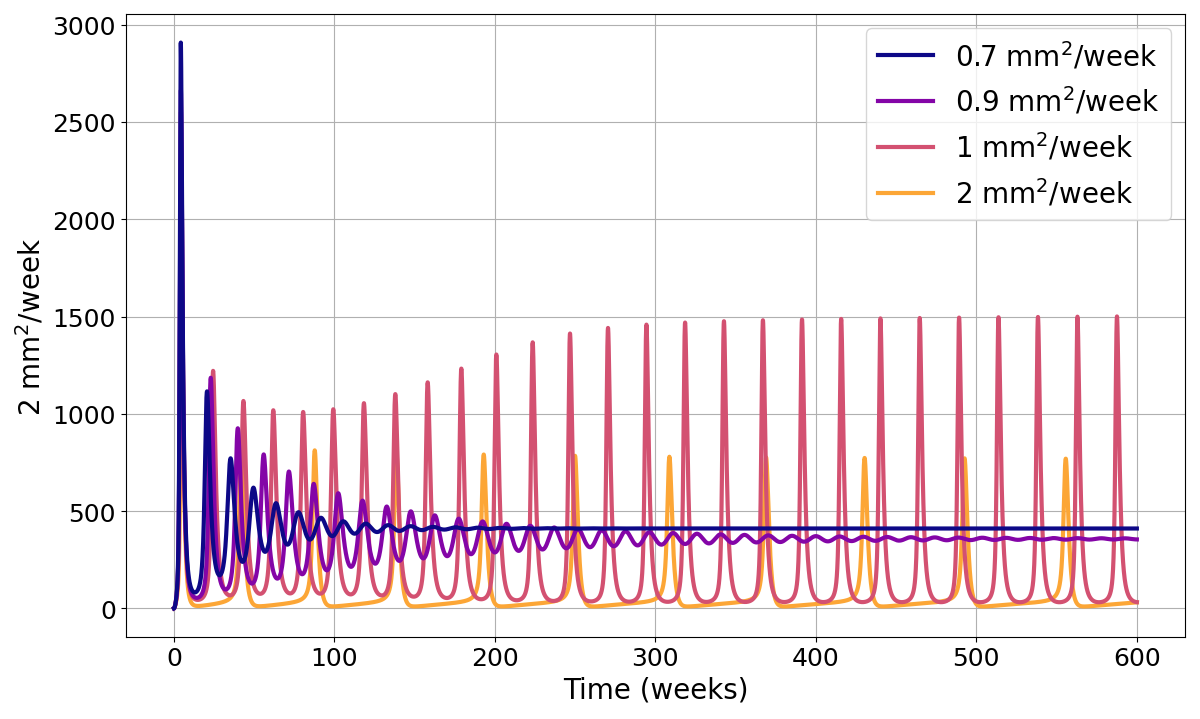}
        \caption{Evolution of the helper lymphocyte concentration in the second cell of a row of $7$ cells over a period of $600$ weeks for four different values of the diffusion coefficient: $\mathscr{D} = 0.7$ \si{mm}$^2$/\si{week}, $\mathscr{D} = 0.9$ \si{mm}$^2$/\si{week}, $\mathscr{D} = 1.0$ \si{mm}$^2$/\si{week} and $\mathscr{D} = 2.0$ \si{mm}$^2$/\si{week}. }
        \label{fig:7_diff_periodic}
\end{figure}

In this situation, the boundary effects alter the behavior of the system, as can be seen in Figure \ref{fig:7_diff_periodic}. Compared to the case of the long row, it is evident that convergence is reached at significantly different times depending on the value of the diffusion coefficient ($\mathscr{D} = 0.7$ \si{mm}$^2$/\si{week} and $\mathscr{D} = 0.9$ \si{mm}$^2$/\si{week}), and the system even seems to exhibit a ``periodic'' trend for higher diffusion coefficients ($\mathscr{D} = 1.0 $ \si{mm}$^2$/\si{week} and $\mathscr{D} = 2.0 $ \si{mm}$^2$/\si{week}).

The value $\mathscr{D} = 1.0$ \si{mm}$^2$/\si{week} is particularly significant, as it marks the boundary (with this level of precision) between a trajectory that converges to a stable fixed point and cyclic behavior. From the figure, we observe that after an initial phase in which the oscillations appear to dampen, even if very slowly, the amplitude then begins to increase again until it stabilizes. When $\mathscr{D} = 2.0$ \si{mm}$^2$/\si{week}, however, we do not observe such evident differences in the values of the peaks following the first one \footnote{There are actually slight differences between the peaks; specifically, it seems that the maximum decreases with each peak. However, these differences are negligible and imperceptible in the image.}, but there is a gradual, albeit slow, increase in the time interval between one peak and the next. The distance between these peaks and their amplitude is dependent on the intensity of diffusion.

The pattern persists for very long times, well beyond those of physical interest, but for these values of $\mathscr{D}$, we can practically consider it periodic, whereas for even higher values, this definition is not justified.

\begin{figure}[ht!]
    \centering
    \includegraphics[width=1.0\linewidth]{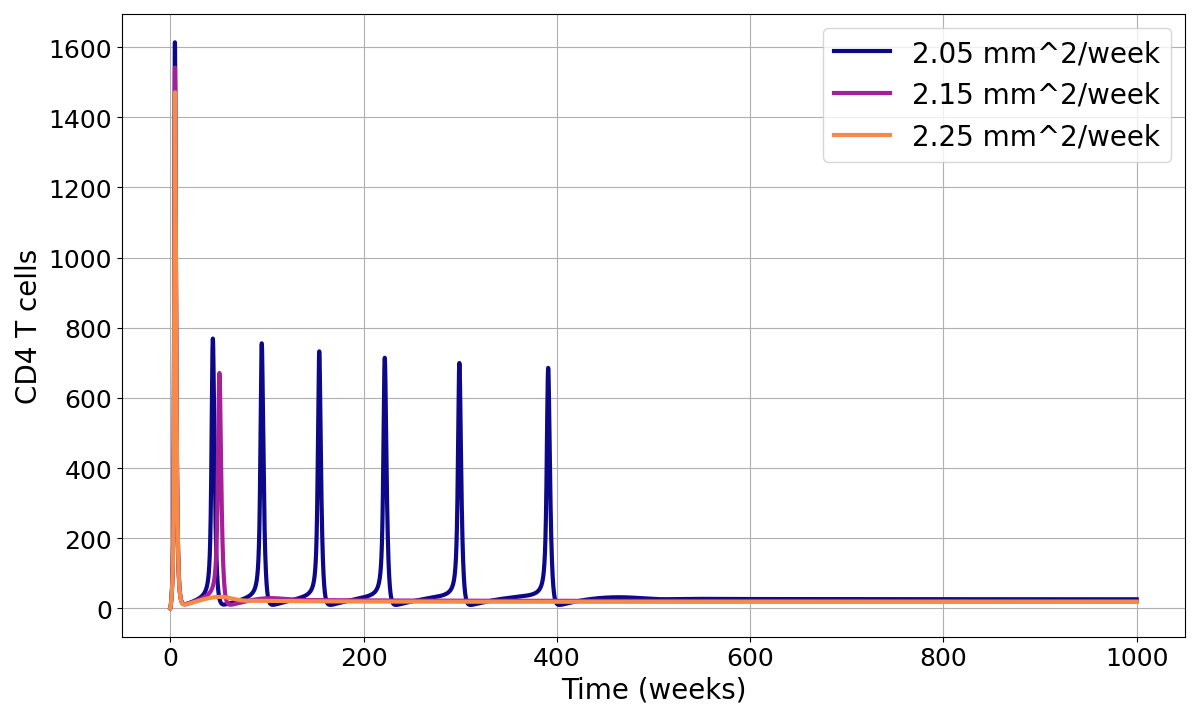}
        \caption{Evolution of the helper lymphocyte concentration (\si{cells/mm}$^3$) in the second cell of a row of $7$ cells over a period of $1000$ weeks for three different values of the diffusion coefficient: $\mathscr{D} = 2.05$ \si{mm}$^2$/\si{week}, $\mathscr{D} = 2.15$ \si{mm}$^2$/\si{week}, and $\mathscr{D} = 2.25$ \si{mm}$^2$/\si{week}. The system eventually converges to a stable fixed point, but oscillations persist for very long times.}
        \label{fig:7_high_diff}
\end{figure}

If we further raise the value of the diffusion coefficient, the tendency to observe increasingly distant peaks becomes extreme, until the pattern eventually fades out before reaching the $1000$-week limit we set for our simulations. 
To rule out the possibility that these oscillations are due to a numerical error, simulations were also performed with different time steps from the usual one. The results remains robust with respect to the choice of $dt$.
We can therefore identify three different types of evolution: damped oscillations that converge to a stable state (for low diffusion), approximately periodic oscillations within physically relevant timescales (for intermediate values of $\mathscr{D}$), and rapid convergence to a stable state, possibly with a limited number of peaks of similar amplitude (for even more intense diffusion). An example of the latter type of evolution is shown in Figure \ref{fig:7_high_diff}, which shows the trend of a single state variable, specifically the concentration of helper lymphocytes, for higher values of the diffusion coefficient. As $\mathscr{D}$ increases, the number of secondary peaks decreases until only a single peak remains. All these trends are qualitatively similar in the other inner cells of the row.

For the cases in which the
muscle is not completely homogeneous at birth, various initializations have been tested, for the three previously defined boundary conditions. 
Our results confirm what was already observed in the uniform case; namely, that diffusion can play under specific boundary and initial conditions, together with particular values of diffusion coefficient. It is worth however noting in detail what happens 
in the presence of a localised injury.

In the model we are studying, muscle damage starts at birth and is represented by the function $\alpha(t)$, which follows a lognormal distribution over time.
Realistically, whether due to wear and tear or, more notably, in the case of trauma, damage is initially localised.
We examine the following scenario:
at birth, the percentage of damaged fibers only in a limited number of cells is greater than zero. Specifically, we consider the simple case where the damaged muscle region (damaged because either $D_0 \neq 0\%$ or $\alpha(t)\neq 0\%$) consists of adjacent cells. 
We impose periodic boundary conditions and simulate the localized damage in a row made of many cells. In this way, damage propagation does not interfere with potential boundary effects that could arise.
Furthermore, due to the cyclic nature of periodic boundary conditions, the position of the damaged tissue within the row is irrelevant.




\begin{figure}[ht!]
\begin{center}
    \begin{minipage}{1.0\textwidth}
        \centering
        \includegraphics[width=1.0\linewidth, trim=0 55 0 0, clip]{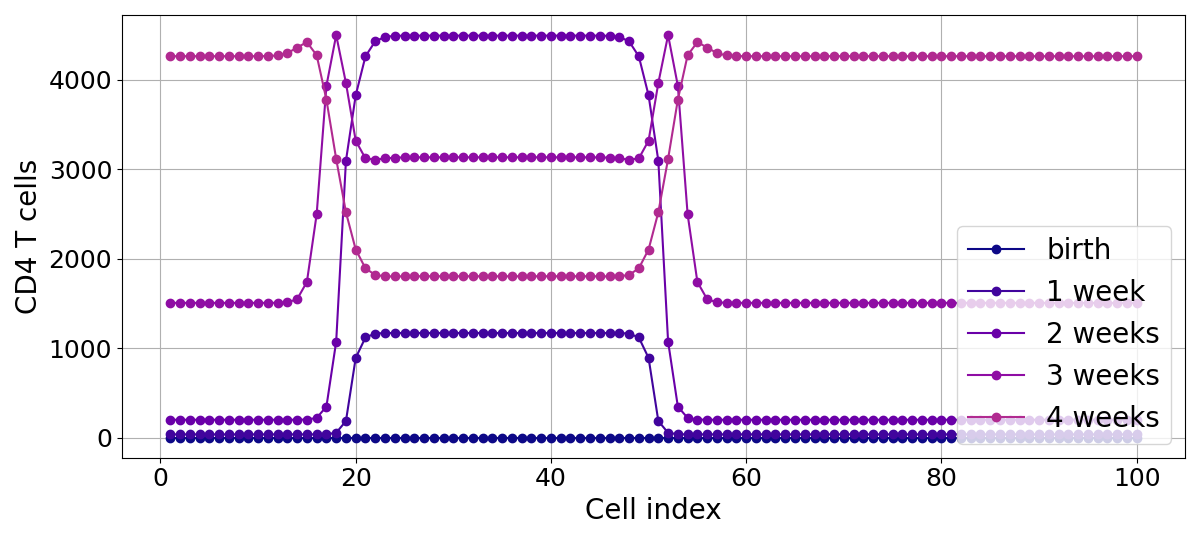}
        \captionsetup{labelsep=none}
    \end{minipage}

    \begin{minipage}{1.0\textwidth}
        \centering
        \includegraphics[width=1.0\linewidth, trim=0 0 0 10, clip]{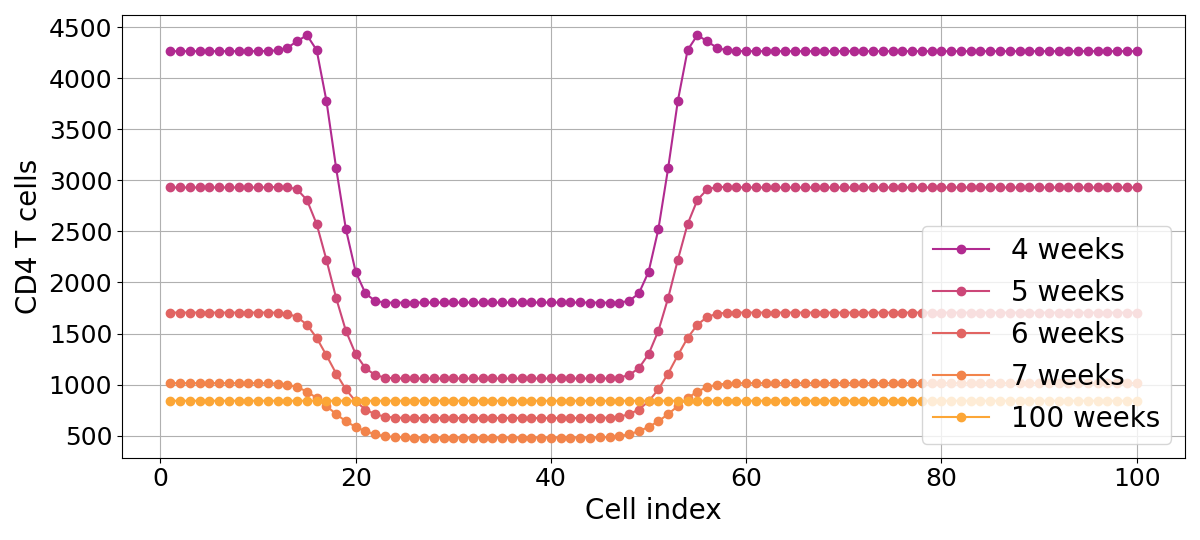}
        \captionsetup{labelsep=none}
    \end{minipage} 
    \caption{Concentration of helper lymphocytes, in \si{cells/mm}$^3$, at various times of the evolution as a function of the cube index (from $1$ to $100$).
    The cubes with a perturbed initial composition range from $20$ to $50$. The boundary conditions are periodic.
    The top plot shows the configurations during the first few weeks, while the lower plot reveals the converge to the steady state.
    The diffusion coefficient is $\mathscr{D}=0.5$ \si{mm}$^2$/\si{week}.
}
    \label{fig:Di_10_01234567100}
\end{center}
\end{figure}


We consider a row of $100$ cubes, where a subset $\mathcal{A}$ is initialized with the default conditions, while for the remaining consecutive cells, subset $\mathcal{B}$, we modify the initialization of the variables $N$, $D$, and $R$. As we choose periodic boundary conditions, the positioning of this subset $\mathcal{B}$ is irrelevant.
As for the initial values $N_0$, $D_0$, and $R_0$, different possibilities are explored for $D_0$, while keeping $R_0 = 0\%$ constant. In reality, tissue damage can already occur in the fetus, and muscle regeneration is a process that takes place even before birth, albeit through different mechanisms. However, we focus here on damage occurring specifically during birth, ensuring that the body, and in particular the immune system, has not yet had the chance to respond to trauma.

The two parts $\mathcal{A}$ and $\mathcal{B}$ are characterized by a different initial configuration. The initial conditions for $H$, $C$, and $M$ remain unchanged. A homogeneous steady state is expected, but how it is realized in time is most relevant for any patient. 
Let us consider a system with $100$ cells, where $D_0^{(i)} = 10\%$ and $N_0^{(i)} = 90\%$ for each index $i$ in the range $[20, 50]$. The diffusion coefficient is $\mathscr{D} = 0.5$ \si{mm}$^2$/\si{week}. The initial localized inhomogeneity gradually fades away, cf.\ Figure \ref{fig:Di_10_01234567100}, that shows the values of a single system variable (the concentration of helper lymphocytes) for each cube in the row at different time points. Each color corresponds to a specific time, allowing us to observe how the response of all the cells in the row evolves over time. At the top, the coordinate is shown during the initial phases of evolution from birth to the fourth week. Below, the values for weeks fourth to seventh are highlighted, followed by those at $100$ weeks, where the system has reached its steady state in each cube.
%
%
We observe that after an initial phase in which the evolution strongly depends on the cube index, the system gradually becomes uniform. 
We also note that, regardless of the time position, a significant portion of the cubes in both intervals $\mathcal{A}$ and
$\mathcal{B}$ remains unaffected by the presence of differently initialized cubes. The number of cells near the initial discontinuity that are influenced by the presence of two different initializations depends on the diffusion coefficient and remains largely unaffected by variations in $\mathcal{B}$.
However, when the width of $\mathcal{B}$ is reduced, the central plateau fails to form, and dips or peaks appear at the corresponding positions. These features become progressively narrower until the limiting case where $\mathcal{B}$ consists of just a single cube.

\section{Two-dimensional model with diffusion}

Consider a square  $N_c \times N_c$ grid of cubes representing a section of tissue. 
The cubes now have neighbors in two dimensions to their right, their left, their top and bottom, through which immune cell transport can occur, as 
given by
the discrete Laplacian in two dimensions:
\begin{align} 
& \mathcal{D}(H_{i,j}) = \frac{D_H}{dx^2} (H_{i+1,j} + H_{i-1,j} + H_{i,j+1} + H_{i,j-1} -4H_{i,j}) \\[0.2cm]
& \mathcal{D}(C_{i,j}) = \frac{D_C}{dx^2} (C_{i+1,j} + C_{i-1,j} + C_{i,j+1} + C_{i,j-1} -4C_{i,j}) \\[0.2cm]
& \mathcal{D}(M_{i,j}) = \frac{D_M}{dx^2} (M_{i+1,j} + M_{i-1,j} + M_{i,j+1} + M_{i,j-1} - 4M_{i,j}),
\end{align}
where $i$ and $j$ are respectively the row and column indices that define the position in the two-dimensional grid.

As in the one dimensional case, three different boundary conditions have been implemented. The simplest initialization is homogeneous across the grid, with fixed concentrations of immune cells at the boundary. 
We now proceed to assess how the results obtained in the one-dimensional case generalize to the system in two dimensions.



We begin with a grid of side length $N_c = 50$.  In the analysis of the $50$-cell row, we assessed the impact of different diffusion coefficients on the system, focusing on the behavior of the state variables in the second cell of the row,  the one closest to the boundary and therefore most affected by diffusion. In the new two-dimensional configuration, this corresponds to studying the evolution of cell $(2, 2)$, or equivalently $(2, 24)$,  $(24, 2)$, or $(24, 24)$, while varying
$\mathscr{D}$. The results obtained are qualitatively similar to those of the one-dimensional case, but the effect of diffusion is more pronounced here, as the cell is close to two boundaries.

\begin{figure}[ht!]
    \centering
    \includegraphics[width=1.0\linewidth, trim=0 0 0 10, clip]{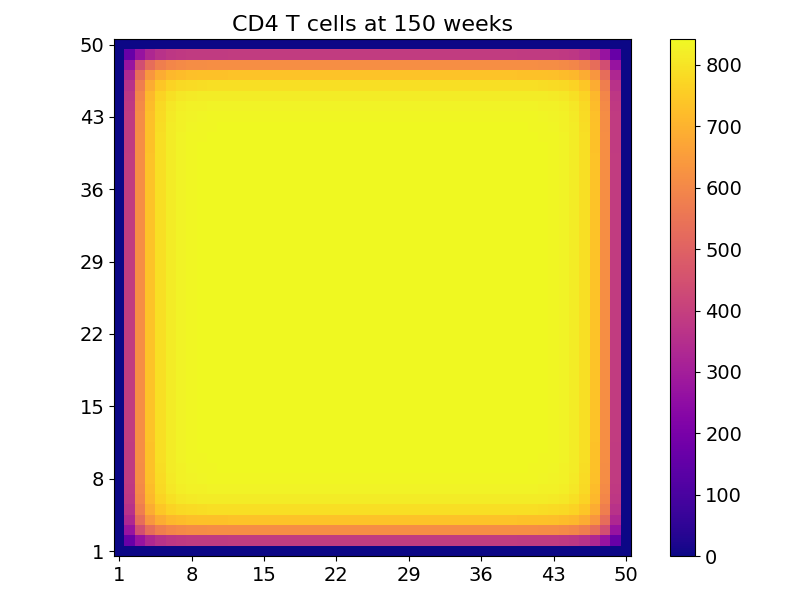}
    \caption{Configuration of the helper lymphocyte concentration across $50 \times 50$ cubes, at steady state ($150$ weeks). The diffusion coefficient is $0.9$ \si{mm}$^2$/\si{week}. The fixed-concentration boundary has a long-term impact limited to the cubes near the edge.}
    \label{fig:2D_50}
\end{figure}

Figure \ref{fig:2D_50} shows the system configuration at long times: $150$ weeks, that is sufficient to reach the steady state. The image shows the heat map of $H$. The chosen diffusion coefficient is $0.9$ \si{mm}$^2$/\si{week}.
Despite this relatively high value, the  one-dimensional phenomenology is confirmed: reaction dominates over diffusion.
Indeed, only the cubes near the fixed-concentration boundary are affected by its presence, while the bulk remains unaffected by this constraint.

\begin{figure}[ht!]
    \centering
    \includegraphics[width=0.95\linewidth, trim=70 0 25 0, clip]{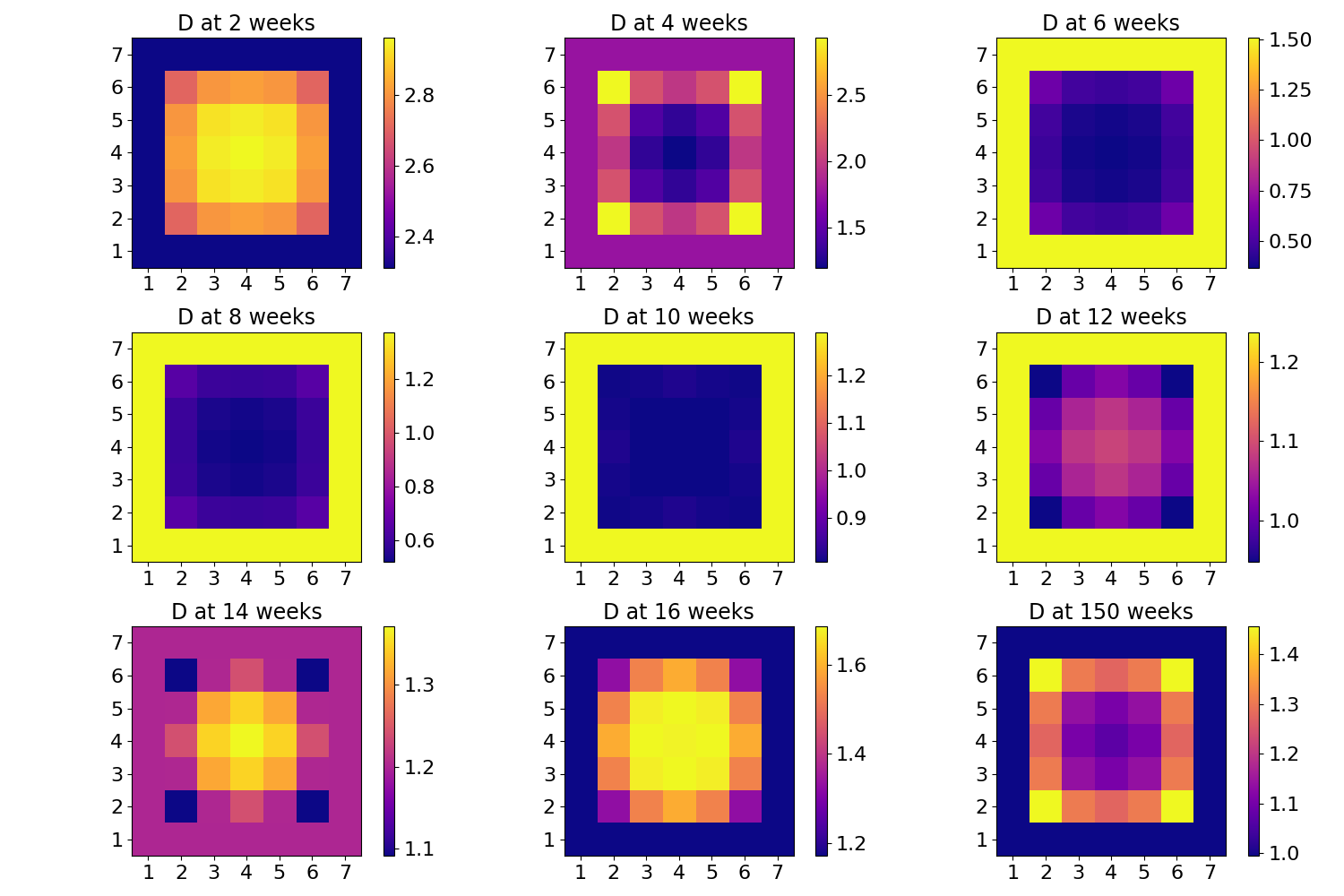}
    \caption{Percentage of damaged fibers in the $7 \times 7$ cubes of the system, at different times. The concentrations are fixed at the boundary, and the diffusion coefficient is $\mathscr{D} = 0.3$ \si{mm}$^2$/\si{week}.}
    \label{fig:D_evol_7cells}
\end{figure}

As in the one-dimensional case, boundary effects become dominant when the system consists of a small number of cells. For instance, consider a $7 \times 7$ grid.
We start with a diffusion coefficient of $0.3$ \si{mm}$^2$/\si{week}. Figure \ref{fig:D_evol_7cells} shows the percentage of damaged fibers in different cubes at various times, ranging from $2$ to $150$ weeks.
This specific variable was chosen for display because its temporal evolution is most evident in the heat map visualization. As expected, even after a long period, the configuration of the grid cells depends on their position: in such a small lattice all the cells are influenced by boundary effects.

Interestingly, the percentage of damaged fibers does not necessarily increase or decrease ``monotonically'' as one moves away from the boundary. For example, in the first image in the top left corner, that is, the snapshot of coordinate $D$ at $2$ weeks, the boundary has the lowest fraction of damaged fibers, which increases as one moves towards the center of the grid. However, in the heat map at $10$ weeks, the opposite occurs. The representation in the bottom left corner (at $14$ weeks) shows yet another case, where cube $(2,2)$ is less damaged than $(1,1)$, but has a lower fraction of damaged fibers than cell $(3,3)$ as well.

In the one-dimensional case, for systems consisting of a few cells with a fixed concentration at the boundary, it was observed that for certain diffusion coefficients, specific patterns emerge in the system's evolution.
To determine whether this phenomenon occurs also in two dimensions, several simulations have been conducted.
Indeed, for certain combinations of $N_c$ and $\mathscr{D}$, the state variables exhibit an oscillatory pattern that does not appear to be damped. 
In the $50 \times 50$ grid example, we observed that the same diffusion coefficient in two dimensions affects the dynamics in cube $(2,2)$ more than it does for cube $2$ in one dimension. 

\begin{figure}[ht!]
    \centering
    \includegraphics[width=1.0\linewidth]{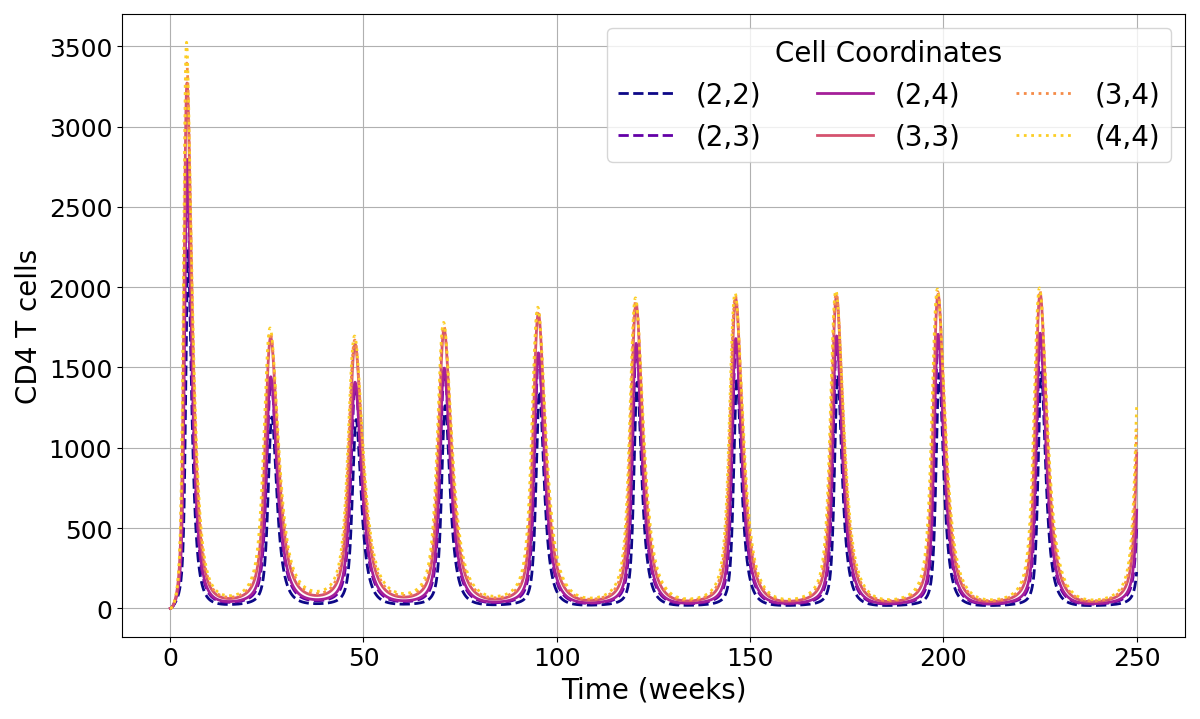}
    \caption{Evolution of helper T cell concentration, in \si{cells/mm}$^3$, in selected central cubes of a $7 \times 7$ grid. The boundaries are maintained at a fixed concentration, and the diffusion coefficient is $0.55$ \si{mm}$^2$/\si{week}.}
    \label{fig:2D_pattern_7cells}
\end{figure}
Consider, for example, a diffusion coefficient of $0.55$ \si{mm}$^2$/\si{week} in a $7 \times 7$ grid.
%
Figure \ref{fig:2D_pattern_7cells} shows all six evolutions for the variable $H$ over $250$ weeks. A trend similar to that of the one-dimensional case is observed: with these boundary conditions and parameters, the variables exhibit nearly periodic oscillations of different amplitude.
In the bulk, the same pattern is observed, which, in the one-dimensional case, would have only emerged for a diffusion coefficient of at least $0.97$ \si{mm}$^2$/\si{week}.
As in the one-dimensional case, we notice that the initial peaks are not perfectly regular, but they begin to stabilise as the evolution progresses. The concentration peaks are higher further away from the boundary.

%% file: Conclusions.tex
\section*{Conclusions and future directions}

Duchenne Muscular Dystrophy (DMD) is a complex condition involving numerous chemical, mechanical, and immunological factors. Through mathematical and numerical modelling, we can now complement experimental studies and enhance our understanding of the processes that characterize it.
In this work, we developed the analysis of a model introduced in Ref.\
\cite{jarrah} concerning the immune response to muscle damage in the \textit{mdx} mouse.
That work provides an idealized description of the role of certain immune cells in the regeneration process. In ours, the thirteen parameters of the model were set according to existing data from the literature.

Our analysis of the model, comprising six ordinary differential equations, was extended beyond the framework of Ref.\ \cite{jarrah}, with particular focus on long-term behavior and the impact of different initial conditions on its evolution. This analysis highlighted a strong dependence of the system on the initial concentration of immune cells. Different values of macrophage and lymphocyte concentrations led to outcomes that were not only quantitatively but also qualitatively distinct from those obtained when using the conditions suggested in the original paper. This indicates a strong connection between the number of immune cells in the tissue and the progression of the disease.

Interestingly, for some sets of initial conditions, the long-term behavior of the system, usually converging to a stable steady state, changes in nature. The values of the state variables begin to oscillate in a periodic-like pattern. This observation aligns with experimental findings that highlight a chronic situation in which regeneration and damage alternate cyclically in the \textit{mdx} mouse. 

Subsequently, we introduced spatial dependence, considering the diffusion of immune cells within the tissue. The one-dimensional case has been studied in detail, and some simulations were also conducted for the two-dimensional version of the model. A close analogy between one and two dimensional cases was found, motivated by the limited impact of diffusion.
Indeed, reaction dominates the evolution in most of the scenarios considered. 
The limited effect of diffusion was further confirmed in the context of localized damage. Although experimental data suggest intense recruitment of immune cells in damaged areas of muscle, diffusion alone does not sufficiently explain this increase in concentration. It is therefore likely that immune cell transport is governed by different mechanisms.
Nonetheless, there are conditions in which diffusion combined with specific boundary effects, leads to a distinctive oscillatory evolution. This phenomenon was first observed in one dimension and persists in the two-dimensional case.


Extensions of the present work could
explore alternative transport mechanisms. For instance, experimental observations indicate that macrophages migrate towards sites of damage via chemotaxis, a transport process driven by specific chemical signals.
Another aspect that remains underexplored is the analysis of the system while accounting for cell crowding. In damaged muscle, the assumption of sparsity, often taken for granted in reaction-diffusion studies, no longer holds. The accumulation of cellular debris and immune cells significantly alters how the system’s components move and interact. A more realistic model would need to incorporate this phenomenon.


\section*{Acknowledgements}
D.C. received government funding managed by the National Research Agency under the France 2030 program, reference ANR-23-IACL-0008.
The support of Italian National Group of
Mathematical Physics (GNFM) of INDAM is acknowledged. L. R.\ and A.A.\ are 
grateful for support from the Italian
Ministry of University and Research (MUR) through
the grant PRIN2022-PNRR project (No. P2022Z7ZAJ)
“A Unitary Mathematical Framework for Modelling
Muscular Dystrophies” (CUP: E53D23018070001).

%% file: appendix1.tex
\section{Damage Modeling and Sensitivity Analysis}
\label{app:A}

Damage is expressed as
\begin{equation}
    \alpha(t) = \frac{h}{t \sigma \sqrt{2 \pi}} \exp\left(- \frac{(\ln{(t)} - m)^2}{2 \sigma^2}\right).
    \label{eq:alpha}
\end{equation}
The functional form of $\alpha(t)$ is designed to model the mechanical damage that arises from a multiplicative degradation process. The basic idea is to consider a random shock $\epsilon_t$ that occurs over time and acts at the level of a single muscle cell. The healthy portion of the system at time $t$, denoted $y_t$, is a fraction of the amount at the previous time step, depending on the shock. As the damage process progresses, the workload and stress on the remaining intact tissue increase progressively as healthy tissue deteriorates. This leads to a feedback loop of accelerated damage. At time $t$, this cumulative effect is given by
\begin{align}
    y_t & = y_{t-1} (1-\epsilon_t) \notag \\
    & = y_0 \prod_{i = 1}^t(1-\epsilon_i).
\end{align}
Taking the logarithm and considering the shocks to be small,
\begin{equation}
    \ln{(y_t)} = \ln{(y_0)} + \sum_{i = 1}^t \ln{(1 - \epsilon_i)} \approx \ln{(y_0)} - \sum_{i = 1}^t \epsilon_i.
\end{equation}
By the central limit theorem, this last sum over $i$ approaches a normal distribution, which implies that $\ln{(y_t)}$ is normally distributed. Since the degradation process determines how quickly $y_t$ drops below the failure threshold, the failure time inherits the same lognormal distribution.

We can rewrite Equation \eqref{eq:alpha} in a more compact and informative form as
\begin{equation}
    \alpha(t) = h \, f_{LN}(t, m, \sigma) \,,
\end{equation}
where
\begin{equation}
    f_{LN}(t, m, \sigma) = \frac{1}{t \sigma \sqrt{2 \pi}} \exp\left(- \frac{(\ln{(t)} - m)^2}{2 \sigma^2}\right)
\end{equation}
is the lognormal probability density function.
Since $f_{LN}(t, m, \sigma)$ integrates to 1, it follows that $h$ is a scale parameter representing the total amount of available damage.
On the other hand, $m$ is the mean of $\log(T)$ when $T$ is a random variable; it therefore shifts the profile on the log-time axis and sets the characteristic time scale of the process (i.e., when damage typically occurs). Finally, $\sigma$ controls the spread in log-time, determining how peaked versus broad the damage profile is. As $\sigma$ increases, damage accumulates more gradually and late-time damage becomes non-negligible (a heavier tail).

\begin{figure}[ht!]
    \centering
    \includegraphics[width=\textwidth]{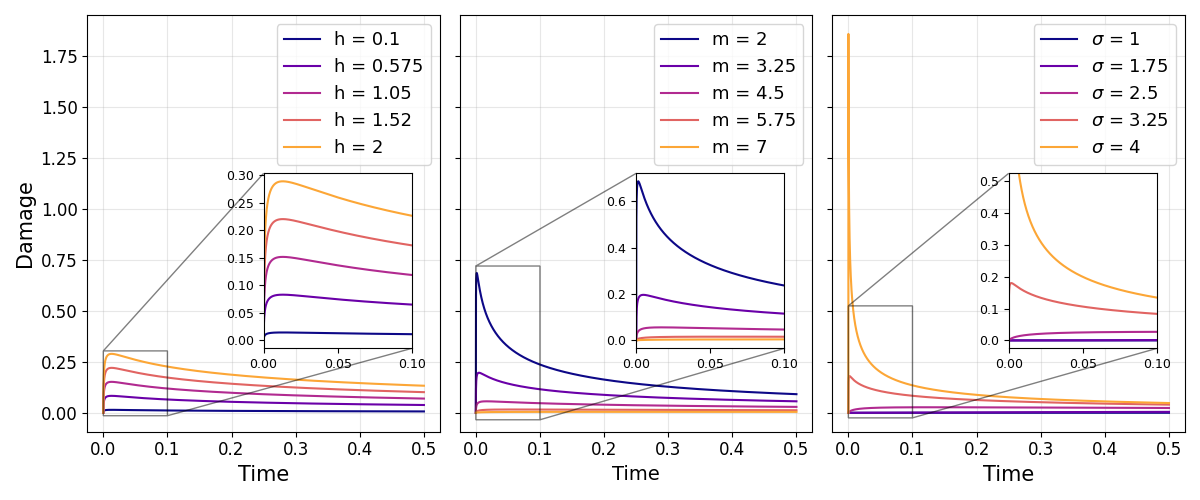}
        \caption{The figure shows how the damage profile depends on the parameters $h, m $ and $\sigma$ (from left to right). The lognormal form in \eqref{eq:alpha} is preserved, and values in a neighborhood of the default ones are explored by varying one parameter at a time.}
        \label{fig:damage1}
\end{figure}

\begin{figure}[ht!]
    \centering
    \includegraphics[width=\textwidth]{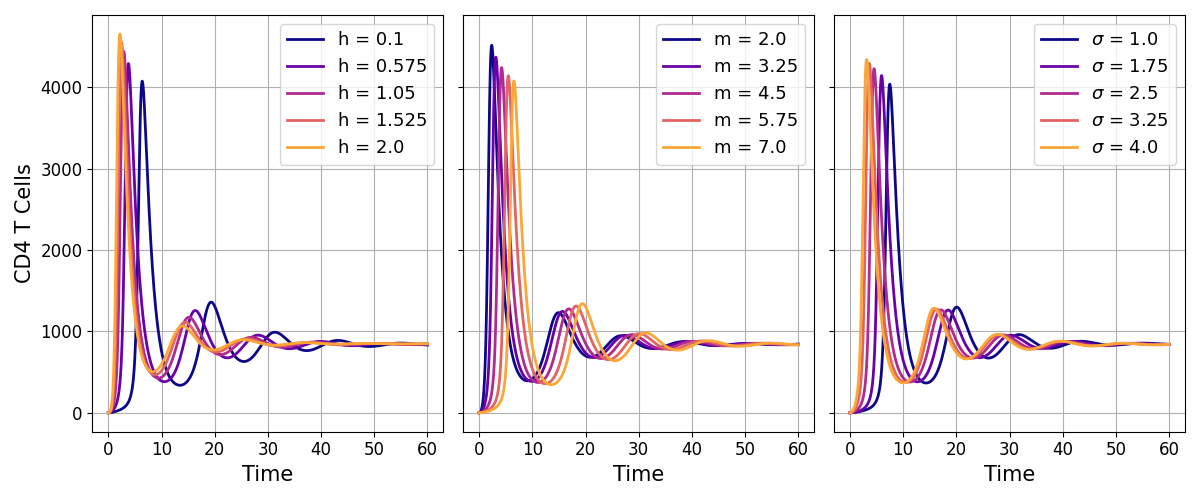}
        \caption{Evolution of the CD4$^+$ T cell concentration for different sets of damage parameters.}
        \label{fig:damage2}
\end{figure}

The default values for the aforementioned parameters are:
\begin{enumerate}
    \item $h_d = 0.511657$
    \item $m_d = 4.22686$
    \item $\sigma_d = 2.92815$.
\end{enumerate}
We now examine how the dynamics change as these parameters vary. Figure \ref{fig:damage1} shows $\alpha(t)$ for different parameter sets.
Consistent with the discussion above, increasing $h$ increases the overall amount of damage (first plot on the left), while increasing $\sigma$ spreads damage over a wider time interval (i.e., the damage profile becomes broader and less peaked). \textit{A priori}, we expect that higher damage leads to a stronger immune response. We test this hypothesis using the same parameter values as in Figure \ref{fig:damage1}, varying one parameter at a time while keeping the others fixed. Figure \ref{fig:damage2} shows the evolution of the helper lymphocyte concentration in the default system as each damage parameter is varied independently.
Over the explored ranges, the trend is qualitatively similar across all parameter sets. Both the heights and the locations of the peaks change, but the response remains damped, with a large initial peak followed by smaller ones. By $60$ weeks, the oscillations are already too small to be visible in the figure and the system is approaching a stationary state, which does not appear to depend strongly on $\alpha_{h,m,\sigma}(t)$ (for each curve, the last recorded value of the plotted concentration lies between $832$ and $850$ \si{cells/mm}$^3$).

\begin{figure}[ht!]
\begin{center}
    \begin{minipage}{1.0\textwidth}
        \centering
        \includegraphics[width=1.0\linewidth, trim=0 27 0 0, clip]{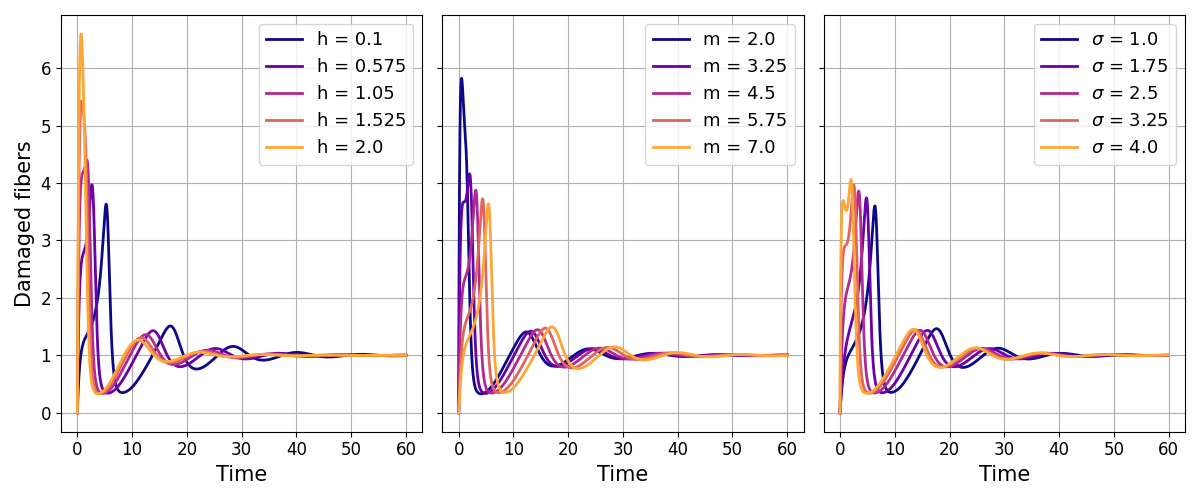}
        \captionsetup{labelsep=none}
    \end{minipage}

    \begin{minipage}{1.0\textwidth}
        \centering
        \includegraphics[width=1.0\linewidth, trim=0 0 0 0, clip]{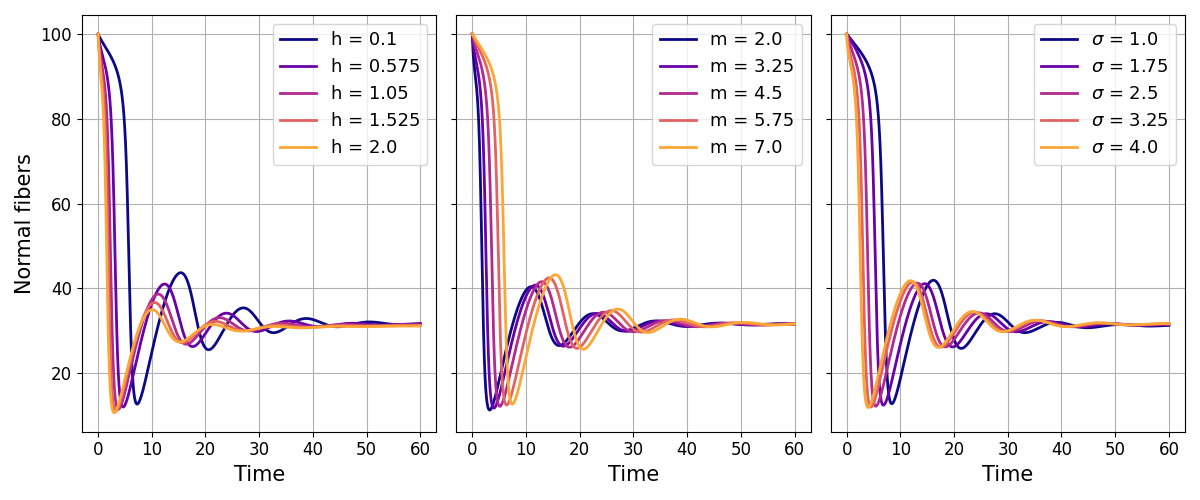}
        \captionsetup{labelsep=none}
    \end{minipage}
    \caption{Evolution of the percentages of damaged (top) and normal (bottom) fibers for different sets of damage parameters.}
    \label{fig:damage3}
\end{center}
\end{figure}

Although changes in immune cell concentrations provide qualitative insight into how differently parameterized damage affects the system, they are not straightforward to interpret quantitatively. It is therefore easier to reason in terms of muscle tissue composition: how do the damage parameters affect the subdivision of muscle fibers into the three classes defined by the model? We are particularly interested in the fraction of damaged fibers under different $\alpha_{h,m,\sigma}(t)$. In simulations with the default $\alpha_{h_d,m_d,\sigma_d}(t)$, this fraction remains fairly low (with a peak around $4\%$) and the healing process does not fail. Damage triggers an immune response that partially regenerates damaged fibers and prevents damage from expanding further. A similar outcome is obtained for the alternative parameter sets, as shown in Figure \ref{fig:damage3}. Focusing on the damaged fibers $D$, the behavior is qualitatively similar across all curves; however, for larger values of $h$ (leftmost panel) the response to $\alpha(t)$ is stronger, and the first peak in damaged fibers exceeds $6\%$ for $h = 2.0$.

Overall, the healing process appears stable with respect to variations in the damage profile. However, intuitively, one might expect regeneration to fail if the fraction of damaged fibers exceeds a critical threshold. To test this, we consider substantially larger values of $h$. Whether $\alpha_{h,m,\sigma}(t)$ remains an appropriate description of damage for such high values of $h$ is not obvious. However, the goal here is not to reproduce a realistic scenario, but rather to assess the system's behavior when a large fraction of fibers is damaged.

\begin{figure}[ht!]
    \centering
    \includegraphics[width=\textwidth]{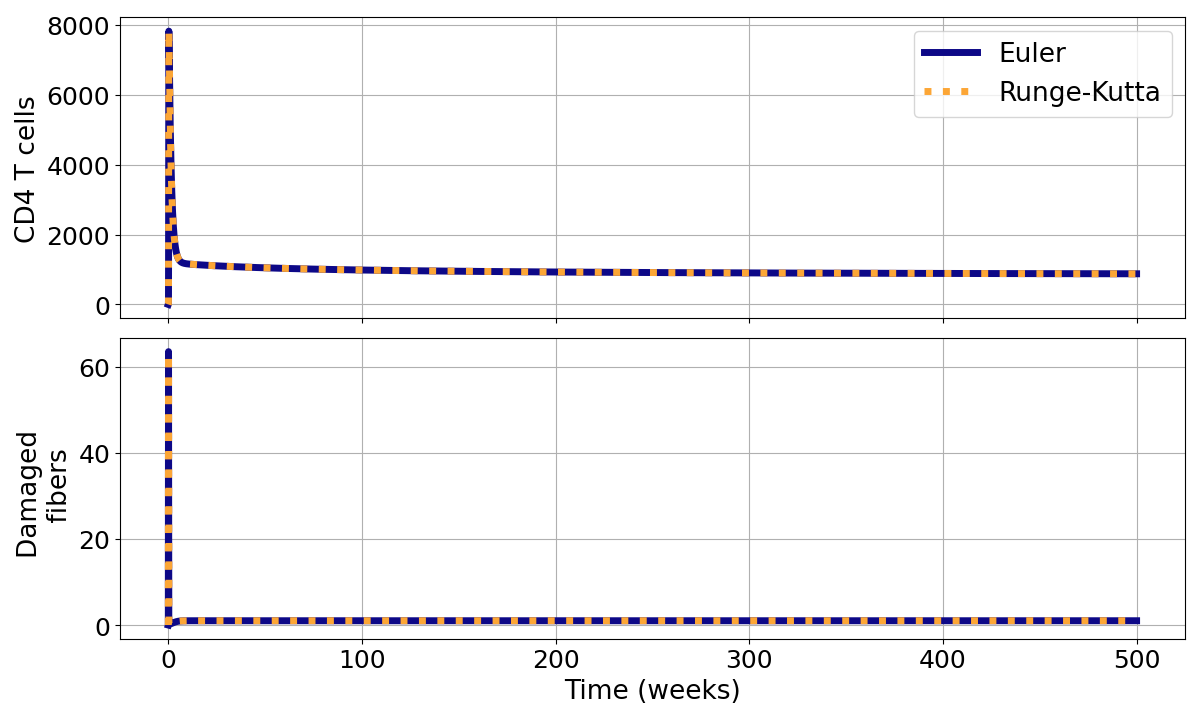}
        \caption{Evolution of the helper lymphocyte concentration (top) and the percentage of damaged fibers (bottom) for $\alpha_{h=100,m_d,\sigma_d}(t)$. For this simulation, $dt=0.001$ weeks.}
        \label{fig:damage5}
\end{figure}
In Figure \ref{fig:damage5}, we show, as an example, the evolution of one immune component (CD4$^+$ T cells concentration) and the percentage of damaged fibers for $h = 100$, $m = m_d$ and $\sigma = \sigma_d$. As expected, the larger damage yields a higher peak in damaged fibers (exceeding $60\%$). Despite the muscle becoming largely damaged, after a sharp peak in the immune response the system stabilizes at lower concentration values. Regenerating fibers account for more than $98\%$ of the total shortly after the peak, and over time the muscle composition approaches that of the long-term default system: at around $1000$ weeks, approximately $31\%$ of the fibers are normal, $68\%$ are regenerating, and only $1\%$ are damaged. The same qualitative behavior is observed for even larger values of $h$, up to $h = 1000$, corresponding to a peak in damaged fibers exceeding $94\%$. It therefore appears that the model predicts an effective immune response across a wide range of $h$.

\section{Positivity of the solutions}
\label{app:pos}
All the state variables of the model~\eqref{eq:H}--\eqref{eq:R} represent
physical quantities (immune-cell concentrations and tissue fractions) and
must therefore remain non-negative throughout the evolution. We verify that
the non-negative orthant
\[
\Omega = \{(H,C,M,N,D,R) : H,C,M,N,D,R \geq 0\}
\]
is positively invariant for the dynamics, i.e.\ any trajectory starting in
$\Omega$ remains in $\Omega$ for all $t > t_0$. Throughout, all parameters
$d_H, d_C, d_M, d_D, k_1, \dots, k_6$ are positive, the resting
concentrations $H_0, C_0, M_0$ are non-negative, and the damage term
$\alpha(t)$ defined in~\eqref{eq:alpha} satisfies $\alpha(t) \geq 0$ for all
$t > 0$.

A standard sufficient condition for positive invariance is that, on each
face of $\Omega$ where a single variable vanishes (while all the others
remain non-negative), the corresponding time derivative points inward, i.e.\
it is non-negative~\cite{haddad2010nonnegative}. The key structural feature of the
model is that every negative contribution to a given equation is
proportional to the variable that is being differentiated, so that all
such terms vanish on the corresponding face. Indeed, evaluating each
equation on its own boundary face gives
\begin{align*}
H = 0 &\;\Rightarrow\; \dot{H} = d_H H_0 + k_1 D M \geq 0, \\
C = 0 &\;\Rightarrow\; \dot{C} = d_C C_0 + k_2 D H \geq 0, \\
M = 0 &\;\Rightarrow\; \dot{M} = d_M M_0 \geq 0, \\
N = 0 &\;\Rightarrow\; \dot{N} = k_4 R \geq 0, \\
D = 0 &\;\Rightarrow\; \dot{D} = k_5 C N + \alpha(t) N \geq 0, \\
R = 0 &\;\Rightarrow\; \dot{R} = k_6 D M + d_D D \geq 0,
\end{align*}
where in each line the remaining variables are taken to be non-negative.
In particular, the depletion terms $-k_5 C N$ and $-\alpha(t) N$
in~\eqref{eq:N}, $-k_6 D M$ and $-d_D D$ in~\eqref{eq:D}, and $-d_H H$,
$-d_C C$, $-d_M M$ in the immune-cell equations, all vanish on the
respective faces precisely because they carry a factor of the vanishing
variable. Since on every face the vector field does not point outward, no
trajectory starting in $\Omega$ can cross its boundary, and the non-negative
orthant is positively invariant.

This guarantees that the model is well posed from a physical standpoint:
the immune-cell concentrations and the muscle-tissue fractions remain
non-negative for all times. Consistently with this analysis, in all our
numerical simulations the state variables were found to remain non-negative
throughout the evolution.

%% file: appendix2.tex
\section{Neumann boundary conditions}
\label{app:B}
The boundary conditions simulated and discussed in this thesis are essentially three: periodic conditions, fixed concentrations at the boundary, and zero flux at the boundary. The latter represents a specific case of a Neumann boundary condition: for a generic concentration $u(x, t)$, one imposes that the flux at the boundary, in the 1D situation, be
\begin{equation*}
    J = D \frac{\partial u}{\partial x} \big|_{boundary} = 0 \, .
\end{equation*}
A Neumann boundary condition specifies the slope of the solution at the boundary. In a biological system, however, situations in which it makes more sense to consider $\partial_x u \neq 0$ are both plausible and indeed rather common. For example, at an interface, the vascular system may act both as a supplier (positive flux) and as a sink or removal mechanism (negative flux). We can therefore use more general Neumann boundary conditions to impose a prescribed influx or efflux at the boundaries.

As a further set of case studies, we considered non-homogeneous Neumann conditions of the type described above, prescribing a nonzero spatial derivative $g_L$ and $g_R$ at the two ends of the domain. These were implemented by correcting the one-sided boundary Laplacian with the imposed gradient contribution, so that in the special case $g_L = g_R = 0$ the scheme reduces to the standard zero-flux (homogeneous Neumann) boundary discretization, in agreement with the original implementation.

\subsubsection*{Simulations}

Depending on the value of the concentration gradient, the flux at the boundary will be stronger or weaker; we denote it by $g_u$, depending on the concentration $u$ to which it refers. This parameter defines the sensitivity of our system to the transport of immune cells from the exterior.

Let us consider the default initial conditions. As already observed, the evolution of the three immune cell populations takes place on very different scales. Consequently, the same value of $g_u$ does not represent a perturbation of the same relative magnitude for all three variables. For this reason, different numerical values were used in the simulations reported here for the different immune cell types. For instance, we tested higher values for macrophages, which reach much larger concentration peaks. More specifically, the gradient values were chosen to span a range from weak to strong boundary forcing.

The imposed Neumann gradients were selected on the basis of the observed dynamic ranges of the diffusing variables in preliminary simulations. One can first identify a characteristic scale of variation,
\begin{equation*}
    \frac{u_{char}}{L},
\end{equation*}
where $L$ is the length of the domain and $u_{char}$ represents a characteristic concentration, and then define the gradient as a fraction of this quantity, either larger or smaller depending on the strength of the perturbation under consideration.

Several assumptions can be made regarding the value assigned to $u_{char}$. Defining $u_{char}$ on the basis of the initial conditions does not capture the subsequent dynamics and suggests gradient values that were found in the simulations to have a negligible effect on the evolution. We therefore chose $u_{char}$ to be of the order of the peak concentration reached by the corresponding immune cell population. This is an empirical and non-rigorous choice, but it makes it possible to account for the different orders of magnitude reached by the various concentrations over time. With this definition, we then explore gradient values that are large enough to be comparable to the reaction terms, but not so large as to dominate the dynamics completely.

Here we restrict ourselves to particularly simple situations in which the gradient has the same magnitude at both boundaries. Moreover, since the concentrations span very different scales, it is sensible to test gradients separately for each variable; that is, in each simulation we set the gradient to zero for two immune cell populations and vary it only for the third. The results reported here were obtained using a diffusion coefficient of 0.5 \si{mm}$^2$/\si{week}. Since boundary forcing enters roughly as $D g_u \, / \, dx$, larger diffusion coefficients amplify the effects described here, whereas smaller diffusion coefficients weaken them. The system under consideration consists of 11 cells.

We study two cases: one in which the effect of the environment is the same at both boundaries, either with two sources or with two sinks at the ends, and one in which the forcing is asymmetric, so that a preferred flow direction can be identified, with inflow on one side and outflow on the other.

\begin{figure}[ht!]
    \centering
    \includegraphics[width=\textwidth]{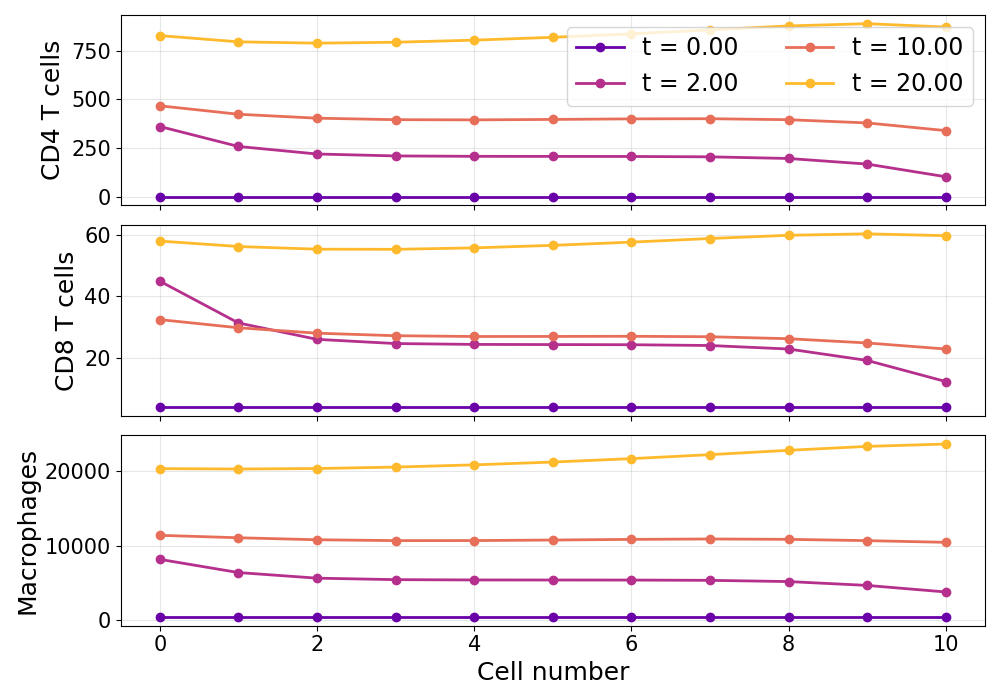}
        \caption{Time series for $|g_H| = 100$. Intuitively, opposite Neumann gradients should create an early imbalance between the two ends of the domain, with lower concentrations near the boundary acting as a sink and higher concentrations near the boundary acting as a source. This behavior is indeed observed in the simulation.}
        \label{fig:neumannH}
\end{figure}

As already mentioned, because the various immune cell populations evolve on different scales, the exploratory ranges over which the effect of the boundary gradient was investigated were as follows:
\begin{align*}
|g_H| & = [0.5, \ldots, 10 ] \\
|g_C| & = [1, \ldots, 100 ] \\
|g_M| & = [50, \ldots, 2000].
\end{align*}

\begin{figure}[ht!]
    \centering
    \includegraphics[width=\textwidth]{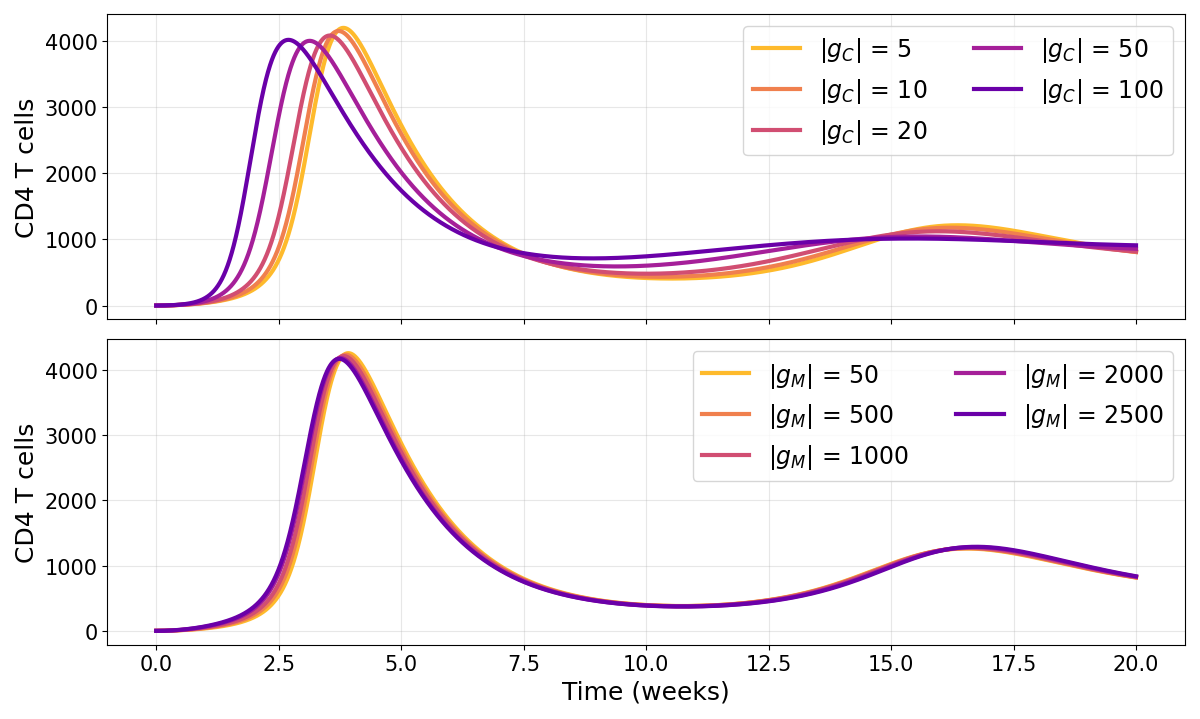}
        \caption{Evolution in the cell at the far right end of the row for a row of 11 cells. At the top, the nonzero boundary gradient is $g_C$, while at the bottom it is $g_M$. Qualitatively, the evolution remains similar regardless of the selected boundary condition.}
        \label{fig:neumannCM_sym}
\end{figure}

However, when the variable under consideration is $H$, these values turn out to be too small to produce an appreciable effect on the evolution of the system. The behavior of the state variables remains very similar across all cells in the row. We therefore conclude that the system is robust with respect to boundary perturbations on $H$ within this explored parameter range. Higher values of $g_H$ were therefore considered in order to identify the regime in which the boundary flux becomes comparable to the intrinsic dynamics of the model. Figure \ref{fig:neumannH} shows an example of the spatial concentration profile of the immune cells in the case of a source of helper lymphocytes (left) and of a sink (right).

\begin{figure}[ht!]
    \centering
    \includegraphics[width=\textwidth]{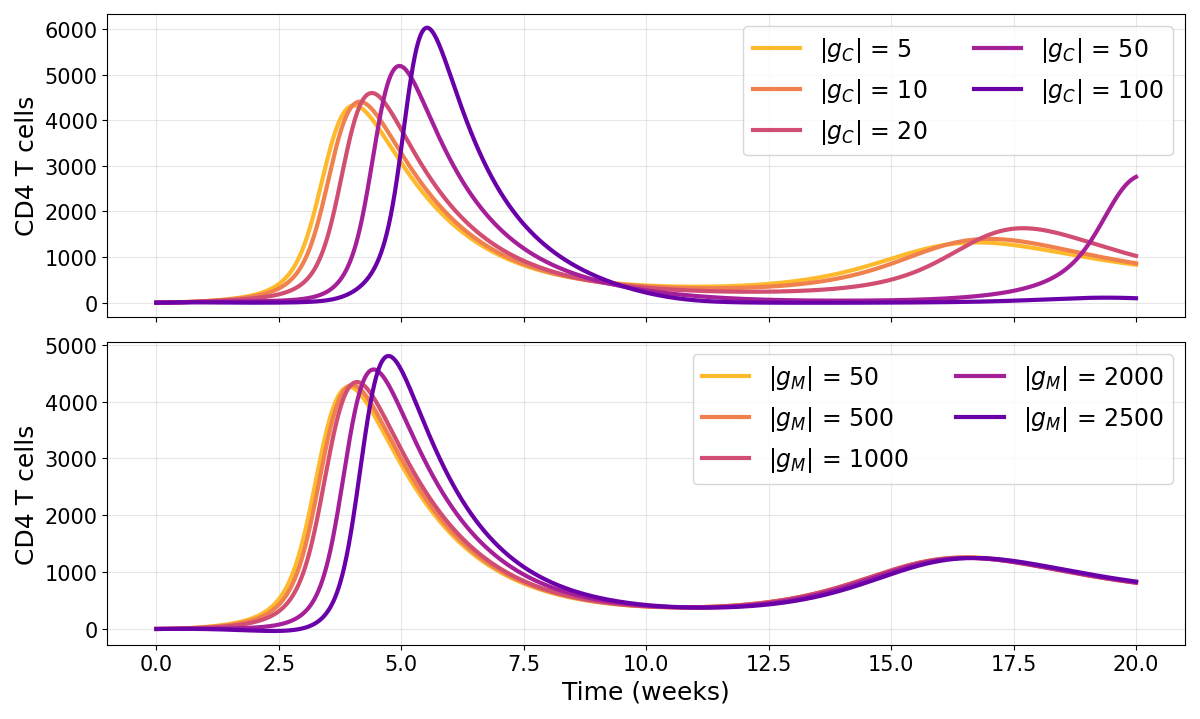}
        \caption{Simulations similar to those in Figure \ref{fig:neumannCM_sym}, but for the case of a source of immune cells at the far left end of the row and a sink at the right end. Different parameter values lead to more pronounced variations in the evolution than in the symmetric case.}
        \label{fig:neumannCM_asym}
\end{figure}

Within the ranges considered for nonzero boundary fluxes associated with the concentrations of cytotoxic lymphocytes and macrophages, differences in the evolution are observed depending on the value of the boundary gradient. In this case as well, however, the qualitative behavior of the system does not change, although, for the same parameter values, the asymmetric case is more sensitive than the symmetric one. It should nevertheless be noted that, under the default initial conditions, some of the gradient values may be too large to be considered plausible. Figure \ref{fig:neumannCM_sym} shows an example of the evolution in the case of two sources at the ends, while Figure \ref{fig:neumannCM_asym} refers to the case of inflow on the left and outflow on the right.

